\documentclass[11pt]{amsart}

\usepackage[francais,english]{babel}
\usepackage{ucs}
\usepackage[utf8]{inputenc}

 \usepackage{hyperref}
 \usepackage{graphicx}

\def\bt{\bf t}

\usepackage{amssymb,epsfig}
\usepackage{bbm}
\usepackage{eufrak}

\newcommand{\R}{\mathbb{R}}

\newcommand{\Z}{\mathbb{Z}}
\renewcommand{\P}{\mathbb{P}}
\newcommand{\Q}{\mathbb{Q}}
\newcommand{\E}{\mathbb{E}}

\newcommand{\N}{\mathbb{N}}

\newcommand{\Tr}{\mathbb{T}}

\newcommand{\B}{\mathbb{B}}
\renewcommand{\Z}{\mathbb{Z}}

\renewcommand{\S}{\mathbb{S}}
\newcommand{\C}{\mathcal{C}}
\newcommand{\Ad}{{\rm A}}
\newcommand{\Bd}{{\rm B}}

\def\build#1_#2^#3{\mathrel{
\mathop{\kern 0pt#1}\limits_{#2}^{#3}}}

\def\cq{$\hfill \square$~\\}
\def \un{\underline}
\def\ind{{\mathbbm{1}}_}

\def\vep{\varepsilon}

\def\lm{\lambda}

\def\ka{\kappa}

\def\rad{{\mathcal R}}
\def\I{{\mathcal I}}
\def\dd{\partial}

\def\m{{\mathcal M}}

\def\f{{\mathcal F}}

\def\u{{\mathcal U}}

\def\t{{\mathcal T}}
\def\v{{\mathcal V}}

\def\W{{\mathcal W}}

\def\d{{\rm d}}
\def\ed{{\rm e}}

\def\z{{\bf z}}
\def\w{{\bf w}}

\def\yg{{\bf y}}
\def\l{{{\mbox{\boldmath$\ell$}}}}
\def\e{{\bf e}}
\def\r{{\bf r}}
\def\b{{\bf b}}
\def\q{{\bf q}}
\def\ag{{\bf a}}
\def\zero{{\bf 0}}
\def\tr{{\bf t}}

\def\map{{\bf m}}

\def\PP{\hbox{\bf P}}

\def\EE{\hbox{\bf E}}

\def\YY{\hbox{\bf Y}}

\def\OO{\mbox{\boldmath$\Omega$}}

\def\mug{\mbox{\boldmath$\mu$}}
\def\zetag{\mbox{\boldmath$\zeta$}}
\def\nug{\mbox{\boldmath$\nu$}}
\def\zetap{{{\mbox{\boldmath$\scriptstyle{\zeta}$}}}}
\def\mup{{{\mbox{\boldmath$\scriptstyle{\mu}$}}}}
\def\nup{{{\mbox{\boldmath$\scriptstyle{\nu}$}}}}
\def\lp{{{\mbox{\boldmath$\scriptstyle{\ell}$}}}}
\def\llbr{[\hspace{-.15em} [ }
\def\rrbr{ ] \hspace{-.15em}]}

\def\be{\begin{equation}}
\def\ee{\end{equation}}
\def\ba{\begin{eqnarray*}}
\def\ea{\end{eqnarray*}}
\def\ov{\overline}
\def\un{\underline}
\def\wh{\widehat}
\def\wt{\widetilde}

\def\la{\longrightarrow}

\def\da{\downarrow}

\def\deg{{\rm deg}}
\def\root{{\varnothing}}
\def\ovla{\overleftarrow}
\def\ovlra{\overleftrightarrow}

\def\proof{\noindent{\bf Proof}\hspace{0.07cm}:\hspace{0.15cm}}

\newtheorem{theorem}{Theorem}[section]
\newtheorem{lemma}[theorem]{Lemma}
\newtheorem{proposition}[theorem]{Proposition}
\newtheorem{corollary}[theorem]{Corollary}

\title{Radius and profile of random planar maps with faces of
  arbitrary degrees}

\author{Grégory Miermont
\and Mathilde Weill
}
\thanks{G.M. --- CNRS, Université de Paris-Sud, Bât. 425, 91405 Orsay
  Cedex, France. \\
 \texttt{gregory.miermont@math.u-psud.fr}} 
\thanks{M.W. --- DMA, \'Ecole normale sup\'erieure, 45 rue d'Ulm, 75005 Paris, France.\\
\texttt{weill@dma.ens.fr}}

\begin{document}

\maketitle

\begin{abstract}
We prove some asymptotic results for the radius and the profile of large
random rooted planar maps with faces of arbitrary degrees. Using a bijection due
to Bouttier, Di Francesco \& Guitter between rooted planar maps and certain
four-type trees with positive labels, we derive our results from a
conditional limit theorem for four-type spatial Galton-Watson trees.
\end{abstract}

\section{Introduction}\label{sec:introduction}

This paper is devoted to the proof of limit theorems for random planar
maps with no constraint on the degree of faces. This work is a natural
sequel to the papers \cite{ChS,LG,MM,We,Mi2}, which dealt with such
limit theorems with an increasing level of generality, starting from
the case of planar quadrangulations and moving to invariance
principles for the radius and the profile of bipartite, then general,
Boltzmann-distributed random planar maps.

Our main goal is to obtain invariance principles for certain
functionals of planar maps with no constraint on the face degrees, of
the same kind as those obtained in \cite{Mi2}. However, while this
work focused on rooted and pointed planar maps, with distances
measured from the distinguished vertex, we focus on maps that are only
rooted and measure distances from the root edge. Similar ``rooted''
results where obtained in \cite{We} building on the ``rooted-pointed''
results of \cite{MM}. 

The basic tools we rely on --- the Bouttier-Di Francesco-Guitter
bijection \cite{BdFG} and methods derived from Le Gall's work
\cite{LG} --- are quite close to those of \cite{We}. However, there
are some notable differences which make the study more intricate. One
of the key differences lies in a change in a re-rooting lemma for
discrete trees, which is considerably more delicate in the present
setting where multiple types are allowed (see Section
\ref{sec:rero-spat-trees}).

Our approach in this paper will be to focus essentially on these
differences, while the parts which can be derived {\it mutatis
  mutandis} from \cite{LG,We} will be omitted.

\section{Preliminaries}\label{sec:preliminaries}

\subsection{Boltzmann laws on planar maps}\label{sec:boltzm-laws-plan}

A {\em planar map} is a proper embedding, without edge crossings, of a
connected graph in the 2-dimensional sphere $\S^2$. Loops and multiple
edges are allowed. The set of vertices will always be equipped with
the graph distance~: if $a$ and $a'$ are two vertices, $d(a,a')$ is
the minimal number of edges on a path from $a$ to $a'$. If $\map$ is a
planar map, we write $\f_\map$ for the set of its faces, and $\v_\map$
for the set of its vertices.

A {\em rooted} planar map is a pair $(\map,\vec{e}\,)$ where $\map$ is
a planar map and $\vec{e}$ is a distinguished oriented edge. The
origin of $\vec{e}$ is called the root vertex. A {\em rooted pointed}
planar map is a triple $(\map,\tau,\vec{e}\,)$ where
$(\map,\vec{e}\,)$ is a rooted planar map and $\tau$ is a
distinguished vertex. We assume that the {\em vertex map}, which is
denoted by $\dag$, is a rooted pointed planar map.

Two rooted maps (resp.~two rooted pointed maps) are identified if
there exists an orientation-preserving homeomorphism of the sphere
that sends the first map to the second one and preserves the root edge
(resp.~the distinguished vertex and the root edge). Let us denote by
$\m_r$ (resp.~$\m_{r,p}$) the set of all rooted maps (resp.~the set of
all rooted pointed maps) up to the preceding identification. In what
follows, we will focus on the subset $\m_{r,p}^+$ of $\m_{r,p}$
defined by~:
$$\m_{r,p}^+=\left\{(\map,\tau,\vec{e}\,)\in\m_{r,p}:d(\tau,e_+)=d(\tau,e_-)+1\right\}\cup\{\dag\}.$$

Let us recall some definitions that can be found in \cite{Mi2}. Let
$\q=(q_i,i\geq1)$ be a sequence of nonnegative weights such that
$q_{2\ka+1}>0$ for at least one $\ka\geq1$. For any planar map
$\map\neq\dag$, we define $W_\q(\map)$ by
$$W_\q(\map)=\prod_{f\in\f_\map}q_{\deg(f)},$$
where we have written ${\rm deg}(f)$ for the degree of the face $f$,
and we set $W_\q(\dag)=1$. We require $\q$ to be admissible that is
$$Z_\q=\sum_{\map\in\m_{r,p}}W_\q(\map)<\infty.$$
Set also
$$Z_\q^+=\sum_{\map\in\m^+_{r,p}}W_\q(\map).$$

For $k,k'\geq0$ we set $N_\bullet(k,k')={2k+k'+1 \choose k+1}$ and
$N_\diamondsuit(k,k')={2k+k'\choose k}$. For every weight sequence we
define 
\ba
f_\q^\bullet(x,y)&=&\sum_{k,k'\geq0}x^ky^{k'}N_\bullet(k,k'){k+k'\choose k}q_{2+2k+k'},\;\;x,y\geq0\\
f_\q^\diamondsuit(x,y)&=&\sum_{k,k'\geq0}x^ky^{k'}N_\diamondsuit(k,k'){k+k'\choose k}q_{1+2k+k'},\;\;x,y\geq0.  
\ea 
From Proposition 1 in \cite{Mi2}, a
sequence $\q$ is admissible if and only if the system 
\ba
\frac{z^+-1}{z^+}&=&f_\q^\bullet(z^+,z^\diamondsuit)\\
z^\diamondsuit&=&f_\q^\diamondsuit(z^+,z^\diamondsuit), \ea has a
solution $(z^+,z^\diamondsuit)\in(0,+\infty)^2$ for which the matrix
${\sf M}_\q(z^+,z^\diamondsuit)$ defined by
$${\sf M}_\q(z^+,z^\diamondsuit)=\left(\begin{array}{ccc} 0 & 0 & z^+-1\\
    \frac{z^+}{z^\diamondsuit}\,\dd_xf_\q^\diamondsuit(z^+,z^\diamondsuit)&\dd_yf_\q^\diamondsuit(z^+,z^\diamondsuit)&0 \\
    \frac{(z^+)^2}{z^+-1}\dd_xf_\q^\bullet(z^+,z^\diamondsuit)&\frac{z^+ z^{\diamondsuit}}{z^+-1}\dd_yf_\q^\bullet(z^+,z^\diamondsuit)&0\end{array}\right)$$
has a spectral radius $\varrho\leq1$. Furthermore this solution is unique and 
\ba
z^{+}&=&Z_\q^+,\\
z^\diamondsuit&=&Z_\q^\diamondsuit,
\ea
where $(Z_\q^\diamondsuit)^2=Z_\q-2Z_\q^++1$. An admissible weight
sequence $\q$ is said to be {\em critical} if the matrix ${\sf
  M}_\q(Z_\q^+,Z_\q^\diamondsuit)$ has a spectral radius
$\varrho=1$. An admissible weight sequence $\q$ is said to be {\em
  regular critical} if $\q$ is critical and if
$f_\q^\bullet(Z_\q^++\vep,Z_\q^\diamondsuit+\vep)<\infty$ for some 
$\vep>0$. 

Let $\q$ be a regular critical weight sequence. We define the
Boltzmann distribution $\B_\q^+$ on the set $\m_{r,p}^+$ by
$$\B_\q^+(\{\map\})=\frac{W_\q(\map)}{Z_\q^+}.$$
Let us now define $Z^{(r)}_\q$ by 
$$Z^{(r)}_\q=\sum_{\map\in\m_r}W_\q(\map).$$
Note that the sum is over the set $\m_r$ of all rooted planar maps. From the fact that $Z_\q<\infty$ it easily follows that $Z^{(r)}_\q<\infty$. We then define the Boltzmann distribution $\B^r_\q$ on the set $\m_r$ by
$$\B^r_\q(\{\map\})=\frac{W_\q(\map)}{Z^{(r)}_\q}.$$

\subsection{The Brownian snake and the conditioned Brownian
  snake}\label{serpent}

Let $x\in\R$. The Brownian snake with initial point $x$ is a pair
$(\b,\r^x)$, where $\b=(\b(s),0\leq s\leq1)$ is a normalized Brownian
excursion and $\r^x=(\r^x(s),0\leq s\leq 1)$ is a real-valued process
such that, conditionally given $\b$, $\r^x$ is Gaussian with mean and
covariance given by
\begin{enumerate}
\item[{$\bullet$}] $\EE[\r^x(s)]=x$ for every $s\in[0,1]$,
\item[{$\bullet$}] ${\bf Cov}(\r^x(s),\r^x(s'))=\displaystyle{\inf_{s\leq t\leq s'}\b(t)}$ for every $0\leq s\leq s'\leq 1$.
\end{enumerate}
We know from \cite{Zu} that $\r^x$ admits a continuous
modification. From now on we consider only this modification. In the
terminology of \cite{Zu} $\r^x$ is the terminal point process of the
one-dimensional Brownian snake driven by the normalized Brownian
excursion $\b$ and with initial point $x$.

Write $\PP$ for the probability measure under which the collection
$(\b,\r^x)_{x\in\R}$ is defined. As mentioned in \cite{We}, for every
$x>0$, we have
$$\PP\left(\inf_{s\in[0,1]}\r^x(s)\geq0\right)>0\, .$$
We may then define for every $x>0$ a pair $(\ov{\b}^x,\ov{\r}^x)$
which is distributed as the pair $(\b,\r^x)$ under the conditioning
that $\inf_{s\in[0,1]}\r^x(s)\geq0$.

We equip $C([0,1],\R)^2$ with the norm $\|(f,g)\|=\|f\|_u\vee\|g\|_u$
where $\|f\|_u$ stands for the supremum norm of $f$. The following
theorem is a consequence of Theorem 1.1 in \cite{LGW}.

\begin{theorem}
  There exists a pair $(\ov{\b}^0,\ov{\r}^0)$ such that
  $(\ov{\b}^x,\ov{\r}^x)$ converges in distribution as $x\da0$ towards
  $(\ov{\b}^0,\ov{\r}^0)$.
\end{theorem} 
The pair $(\ov{\b}^0,\ov{\r}^0)$ is the so-called conditioned Brownian
snake with initial point $0$.

Theorem 1.2 in \cite{LGW} provides a useful construction of the
conditioned object $(\ov{\b}^0,\ov{\r}^0)$ from the unconditioned one
$(\b,\r^0)$. In order to present this construction, first recall that
there is a.s.~a unique $s_*$ in $(0,1)$ such that
$$\r^0(s_*)=\inf_{s\in[0,1]}\r^0(s)$$
(see Lemma 16 in \cite{MaMo2} or Proposition 2.5 in \cite{LGW}). For
every $s\in[0,\infty)$, write $\{s\}$ for the fractional part of
$s$. According to Theorem 1.2 in \cite{LGW}, the conditioned snake
$(\ov{\b}^0,\ov{\r}^0)$ may be constructed explicitly as follows~: for
every $s\in[0,1]$, \ba
\ov{\b}^0(s)&=&\b({s_*})+\b(\{s_*+s\})-2\,\inf_{s\wedge\{s_*+s\}\leq
  t\leq
  s\vee\{s_*+s\}}\b(t),\\
\ov{\r}^0(s)&=&\r^0(\{s_*+s\})-\r^0({s_*}).  \ea

\subsection{Statement of the main result}

We first need to introduce some notation. Let $\map\in\m_r$. We denote
by $o$ its root vertex. The radius $\rad_\map$ is the maximal distance
between $o$ and another vertex of $\map$ that is
$$\rad_\map=\max\{d(o,a):a\in\v_\map\}.$$
The profile of $\map$ is the measure $\lm_\map$ on $\{0,1,2,\ldots\}$
defined by
$$\lm_\map(\{k\})=\#\{a\in\v_\map:d(o,a)=k\},\;k\geq0.$$
Note that $\rad_\map$ is the supremum of the support of $\lm_\map$. It
is also convenient to introduce the rescaled profile. If $\map$ has
$n$ vertices, this is the probability measure on $\R_+$ defined by
$$\lm_\map^{(n)}(A)=\frac{\lm_\map(n^{1/4}A)}{n}$$
for any Borel subset $A$ of $\R_+$.

Recall from section \ref{serpent} that $(\b,\r^0)$ denotes the
Brownian snake with initial point $0$.

\begin{theorem}\label{thcartes}
  Let $\q$ be a regular critical weight sequence. There exists a
  scaling constant ${\rm C}_\q$ such that the following results hold.
\begin{enumerate}
\item[(i)] The law of $n^{-1/4}\,\rad_\map$ under
  $\B^r_\q(\cdot\mid\#\v_\map=n)$ converges as $n\to\infty$ to the law
  of the random variable
$${\rm C}_\q\left(\sup_{0\leq s\leq1}\r^0(s)-\inf_{0\leq s\leq1}\r^0(s)\right).$$
\item[(ii)] The law of the random probability measure $\lm_\map^{(n)}$
  under $\B^r_\q(\cdot\mid\#\v_\map=n)$ converges as $n\to\infty$ to
  the law of the random probability measure $\I$ defined by
$$\langle\I,g\rangle=\int_0^1g\left({\rm C}_\q\left(\r^0(t)-\inf_{0\leq s\leq1}\r^0(s)\right)\right)\d t.$$
\item[(iii)] The law of the rescaled distance $n^{-1/4}\,d(o,a)$ where
  $a$ is a vertex chosen uniformly at random among all vertices of
  $\map$, under $\B^r_\q(\cdot\mid\#\v_\map=n)$ converges as
  $n\to\infty$ to the law of the random variable
$${\rm C}_\q\sup_{0\leq s\leq1}\r^0(s).$$
\end{enumerate}
\end{theorem}

Theorem \ref{thcartes} is an analogue to Theorem 2.5 in \cite{We} in
the setting of non-bipartite maps. Beware that in Theorem
\ref{thcartes} maps are conditioned on their number of vertices
whereas in \cite{We} they are conditioned on their number of
faces. However the results stated in Theorem 2.5 in \cite{We} remain
valid by conditioning on the number of vertices (with different
scaling constants). On the other hand, our arguments to prove Theorem
\ref{thcartes} do not lead to the statement of these results by
conditioning maps on their number of faces. A notable exception is the
case of $k$-angulations ($\q=q\delta_k$ for some $k\geq 3$ and
appropriate $q>0$), where an application of Euler's formula shows that
$\#\f_\map=(k/2-1)\#\v_\map+2$, so that the two conditionings are
essentially equivalent and result in a change in the scale factor
${\rm C}_\q$.

Recall that the results of Theorem 2.5 in \cite{We} for the special
case of quadrangulations were obtained by Chassaing \& Schaeffer
\cite{ChS} (see also Theorem 8.2 in \cite{LG}).

Last observe that Theorem \ref{thcartes} is obviously related to Theorem 1 in
\cite{Mi2}. Note however that \cite{Mi2} deals with rooted pointed
maps instead of rooted maps as we do and studies distances from the
distinguished point of the map rather than from the root vertex.

\subsection{Multitype spatial trees}

We start with some formalism for discrete trees. Set
$$\u=\bigcup_{n\geq0}\N^n,$$
where by convention $\N=\{1,2,3,\ldots\}$ and $\N^0=\{\root\}$. An
element of $\,\u$ is a sequence $u=u^1\ldots u^n$, and we set $|u|=n$
so that $|u|$ represents the generation of $u$. In particular
$\vert\varnothing\vert=0$. If $u=u^1\ldots u^n$ and $v=v^1\ldots v^m$
belong to $\,\u$, we write $uv=u^1\ldots u^nv^1\ldots v^m$ for the
concatenation of $u$ and $v$. In particular $\varnothing
u=u\varnothing=u$. If $v$ is of the form $v=uj$ for $u\in\u$ and
$j\in\N$, we say that $v$ is a child of $u$, or that $u$ is the father
of $v$, and we write $u=\check{v}$. More generally if $v$ is of the
form $v=uw$ for $u,w\in\u$, we say that $v$ is a descendant of $u$, or
that $u$ is an ancestor of $v$. The set $\u$ comes with the natural
lexicographical order such that $u\preccurlyeq v$ if either $u$ is an
ancestor of $v$, or if $u=wa$ and $v=wb$ with $a\in\u^\ast$ and
$b\in\u^\ast$ satisfying $a^1<b^1$, where we have set
$\u^\ast=\u\setminus\{\varnothing\}$. We write $u\prec v$ if
$u\preccurlyeq v$ and $u\neq v$.

A plane tree $\tr$ is a finite subset of $\,\u$ such that 
\begin{enumerate}
\item[$\bullet$] $\varnothing\in\tr$,
\item[$\bullet$] $u\in\tr\setminus\{\varnothing\}\Rightarrow\check{u}\in\tr$,
\item[$\bullet$] for every $u\in\tr$ there exists a number
  $c_u(\tr)\geq0$ such that $uj\in\tr\Leftrightarrow1\leq j\leq
  c_u(\tr)$.
\end{enumerate}

Let $\tr$ be a plane tree and let $\xi=\#\tr-1$. The {\em search-depth
  sequence} of $\tr$ is the sequence $u_0,u_1,\ldots,u_{2\xi}$ of
vertices of $\tr$ wich is obtained by induction as follows. First
$u_0=\root$, and then for every $i\in\{0,1,\ldots,2\xi-1\}$, $u_{i+1}$
is either the first child of $u_i$ that has not yet appeared in the
sequence $u_0,u_1,\ldots,u_i$, or the father of $u_i$ if all children
of $u_i$ already appear in the sequence $u_0,u_1,\ldots,u_i$. It is
easy to verify that $u_{2\xi}=\root$ and that all vertices of $\tr$
appear in the sequence $u_0,u_1,\ldots,u_{2\xi}$ (of course some of
them appear more that once). We can now define the {\em contour
  function} of $\tr$. For every $k\in\{0,1,\ldots,2\xi\}$, we let
$C(k)=|u_k|$ denote the generation of the vertex $u_{k}$. We extend
the definition of $C$ to the line interval $[0,2\xi]$ by interpolating
linearly between successive integers. Clearly $\tr$ is uniquely
determined by its contour function $C$.

Let $K\in\N$ and $[K]=\{1,2,\ldots,K\}$. A $K$-type tree is a pair
$(\tr,\e)$ where $\tr$ is a plane tree and $\e:\tr\to[K]$ assigns a
type to each vertex. If $(\tr,\e)$ is a $K$-type tree and if $i\in[K]$
we set
$$\tr^i=\{u\in\tr:\e(u)=i\}.$$
We denote by $T^{(K)}$ the set of all $K$-type trees and we set
$$T^{(K)}_i=\left\{(\tr,\e)\in T^{(K)}:\e(\root)=i\right\}.$$ 
Set
$$\W_K=\bigcup_{n\geq0}[K]^n,$$
with the convention $[K]^0=\{\root\}$. An element of $\W_K$ is a
sequence $\w=(w_1,\ldots,w_n)$ and we set $|\w|=n$. Consider the
natural projection $p:\W\to\Z_+^K$ where
$p(\w)=(p_1(\w),\ldots,p_K(\w))$ and
$$p_i(\w)=\#\{j\in\{1,\ldots,|\w|\}:w_j=i\}.$$
Let $u\in\u$ and let $(\tr,\e)\in T^{(K)}$ such that $u\in\tr$. We
then define $\w_u(\tr)\in\W_K$ by
$$\w_u(\tr)=(\e(u1),\ldots,\e(uc_u(\tr))),$$
and we set $\z_u(\tr)=p(\w_u(\tr))$.

A $K$-type spatial tree is a triple $(\tr,\e,\l)$ where $(\tr,\e)\in
T^{(K)}$ and $\l:\tr\to\R$. If $v$ is a vertex of $\tr$, we say that
$\l_v$ is the {\em label} of $v$. We denote by $\Tr^{(K)}$ the set of
all $K$-type spatial trees and we set
$$\Tr^{(K)}_i=\left\{(\tr,\e,\l)\in\Tr^{(K)}:\e(\root)=i\right\}.$$
If $(\tr,\e,\l)\in\Tr^{(K)}$ we define the {\em spatial contour
  function} of $(\tr,\e,\l)$ as follows. Recall that
$u_0,u_1,\ldots,u_{2\xi}$ denotes the search-depth sequence of
$\tr$. First if $k\in\{0,\ldots,2\xi\}$, we put $V(k)=\l_{u_k}$. We
then complete the definition of $V$ by interpolating linearly between
successive integers.

\subsection{Multitype spatial Galton-Watson trees}

Let $\zetag=(\zeta^{(i)},i\in[K])$ be a family of probability measures
on the set $\W_K$. We associate with $\zetag$ the family
$\mug=(\mu^{(i)},i\in[K])$ of probability measures on the set $\Z_+^K$
in such a way that each $\mu^{(i)}$ is the image measure of
$\zeta^{(i)}$ under the mapping $p$. We make the basic assumption that
$$\max_{i\in[K]}\mu^{(i)}\Bigg(
\Bigg\{\z\in\Z_+^K:\sum_{j=1}^Kz_j\neq1\Bigg\}\Bigg)>0,$$ and we say
that $\zetag$ (or $\mug$) is non-degenerate. If for every $i\in[K]$,
$\w\in\W_K$ and $\z=p(\w)$ we have
$$\zeta^{(i)}(\{\w\})=\frac{\mu^{(i)}(\{\z\})}{\#\left(p^{-1}(\z)\right)},$$
then we say that $\zetag$ is the {\em uniform ordering} of $\mug$.

For every $i,j\in[K]$, let
$$m_{ij}=\sum_{\z\in\Z_+^K}z_j\mu^{(i)}(\{\z\}),$$
be the mean number of type-$j$ children of a type-$i$ individual, and
let ${\sf M}_\mup=(m_{ij})_{1\leq i,j\leq K}$. The matrix ${\sf
  M}_\mup$ is said to be irreducible if for every $i,j\in[K]$ there
exists $n\in\N$ such that $m_{ij}^{(n)}>0$ where we have written
$m_{ij}^{(n)}$ for the $ij$-entry of ${\sf M}_\mup^n$. We say that
$\zetag$ (or $\mug$) is irreducible if ${\sf M}_\mup$ is. Under this
assumption the Perron-Frobenius theorem ensures that ${\sf M}_\mup$
has a real, positive eigenvalue $\varrho$ with maximal modulus. The
distribution $\zetag$ (or $\mug$) is called sub-critical if
$\varrho<1$ and critical if $\varrho=1$.

Assume that $\zetag$ is non-degenerate, irreducible and
(sub-)critical. We denote by $P^{(i)}_\zetap$ the law of a $K$-type
Galton-Watson tree with offspring distribution $\zetag$ and with
ancestor of type $i$, meaning that for every $(\tr,\e)\in T_i^{(K)}$,
$$P^{(i)}_\zetap(\{(\tr,\e)\})=\prod_{u\in\tr}\zeta^{(\e(u))}\left(\w_u(\tr)\right),$$
The fact that this formula defines a probability measure on
$T_i^{(K)}$ is justified in \cite{Mi1}.

Let us now recall from \cite{Mi1} how one can couple $K$-type trees
with a spatial displacement in order to turn them into random elements
of $\Tr^{(K)}$. To this end, consider a family
$\nug=(\nu_{i,\w},i\in[K],\w\in\W_K)$ where $\nu_{i,\w}$ is a
probability measure on $\R^{|\w|}$. If $(\tr,\e)\in T^{(K)}$ and
$x\in\R$, we denote by $R_{\nu,x}((\tr,\e),\d\l)$ the probability
measure on $\R^\tr$ which is characterized as follows. For every $i\in[K]$ and $u\in\tr$ such that $\e(u)=i$, consider
$\YY_u=(Y_{u1},\ldots,Y_{u|\w|})$ (where we have written
$\w_u(\tr)=\w$) a random variable distributed according to
$\nu_{i,\w}$, in such a way that $(\YY_u,u\in\tr)$ is a collection of
independant random variables. We set $L_\root=x$ and for every
$v\in\tr\setminus\{\root\}$,
$$L_v=\sum_{u\in\,\rrbr\root,v\rrbr}Y_u,$$
where $\rrbr\root,v\rrbr$ is the set of all ancestors of $v$ distinct
from the root $\root$. The probability measure
$R_{\nup,x}((\tr,\e),\d\l)$ is then defined as the law of
$(L_v,v\in\tr)$. We finally define for every $x\in\R$ a probability
measure $\P^{(i)}_{\zetap,\nup,x}$ on the set $\Tr_i^{(K)}$ by
setting,
$$\P^{(i)}_{\zetap,\nup,x}(\d\tr\,\d\e\,\d\l)=P^{(i)}_\zetap(\d\tr,\d\e)R_{\nup,x}((\tr,\e),\d\l).$$

\subsection{The Bouttier-Di Francesco-Guitter bijection}

We start with a definition. We consider the set $T_M\subset T_1^{(4)}$
of $4$-type trees in which, for every $(\tr,\e)\in T_M$ and $u\in\tr$,
\begin{enumerate}
\item[{\bf 1.}] if $\e(u)=1$ then $\z_u(\tr)=(0,0,k,0)$ for some $k\geq0$,
\item[{\bf 2.}] if $\e(u)=2$ then $\z_u(\tr)=(0,0,0,1)$,
\item[{\bf 3.}] if $\e(u)\in\{3,4\}$ then $\z_u(\tr)=(k,k',0,0)$ for
  some $k,k'\geq0$.
\end{enumerate}
Let now $\Tr_M\subset\Tr_1^{(4)}$ be the set of $4$-type spatial trees $(\tr,\e,\l)$ such that $(\tr,\e)\in T_M$ and in which, for every $(\tr,\e,\l)\in\Tr_M$ and $u\in\tr$,
\begin{enumerate}
\item[{\bf 4.}] $\l_u\in\Z$,
\item[{\bf 5.}] if $\e(u)\in\{1,2\}$ then $\l_u=\l_{ui}$ for every $i\in\{1,\ldots,c_u(\tr)\}$,
\item[{\bf 6.}] if $\e(u)\in\{3,4\}$ and $c_u(\tr)=k$ then by setting $u0=u(k+1)=\check{u}$ and $x_i=\l_{ui}-\l_{u(i-1)}$ for $1\leq i\leq k+1$, we have
\begin{enumerate}
\item[(a)] if $\e(u(i-1))=1$ then $x_i\in\{-1,0,1,2,\ldots\}$,
\item[(b)] if $\e(u(i-1))=2$ then $x_i\in\{0,1,2,\ldots\}$.
\end{enumerate}
\end{enumerate}
We will be interested in the set
$\ov{\Tr}_M=\{(\tr,\e,\l)\in\Tr_M:\ell_\root=1\;{\rm
  and}\;\l_v\geq1\;{\rm for}\;{\rm all}\;v\in\tr^1\}$. Notice that
condition {\bf 6.} implies that if $(\tr,\e,\l)\in\ov{\Tr}_M$ then
$\l_v\geq0$ for all $v\in\tr$.

We will now describe the Bouttier-Di Francesco-Guitter bijection from
the set $\ov{\Tr}_{M}$ onto $\m_{r}$. This bijection can be found in
\cite{BdFG} in the more general setting of Eulerian maps.

Let $(\tr,\e,\l)\in\ov{\Tr}_{M}$. Recall that $\xi=\#\tr-1$. Let
$u_0,u_1,\ldots,u_{2\xi}$ be the search-depth sequence of $\tr$. It is
immediate to see that $\e(u_k)\in\{1,2\}$ if $k$ is even and that
$\e(u_k)\in\{3,4\}$ if $k$ is odd. We define the sequence
$v_0,v_1,\ldots,v_\xi$ by setting $v_k=u_{2k}$ for every
$k\in\{0,1,\ldots,\xi\}$. Notice that $v_0=v_\xi=\root$.

Suppose that the tree $\tr$ is drawn in the plane and add an extra
vertex $\dd$, not on $\tr$. We associate with $(\tr,\e,\l)$ a planar
map whose set of vertices is
$$\tr^1\cup\{\dd\},$$
and whose edges are obtained by the following device~: for every
$k\in\{0,1,\ldots,\xi-1\}$,
\begin{enumerate}
\item[$\bullet$] if $\e(v_k)=1$ and $\l_{v_k}=1$, or if $\e(v_k)=2$ and $\l_{v_k}=0$, draw an edge between $v_k$ and $\dd$~;
\item[$\bullet$] if $\e(v_k)=1$ and $\l_{v_k}\geq2$, or if $\e(v_k)=2$
  and $\l_{v_k}\geq1$, draw an edge between $v_k$ and the first vertex
  in the sequence $v_{k+1},\ldots,v_{\xi}$ with type $1$ and
  label $\l_{v_k}-\ind{\{\e(v_k)=1\}}$.
\end{enumerate}
Notice that condition {\bf 6.} in the definition of the set
$\ov{\Tr}_M$ entails that
$\l_{v_{k+1}}\geq\l_{v_k}-\ind{\{\e(v_k)=1\}}$ for every
$k\in\{0,1,\ldots,\xi-1\}$, and recall that
$\min\{\l_{v_j}:j\in\{k+1,\ldots,\xi\}\;{\rm and}\;\e(v_j)=1\}=1$. The
preceding properties ensure that whenever $\e(v_k)=1$ and
$\l(v_k)\geq2$ or $\e(v_k)=2$ and $\l(v_k)\geq1$ there is at least one
type-$1$ vertex among $\{v_{k+1},\ldots,v_\xi\}$ with label
$\l_{v_k}-\ind{\{\e(v_k)=1\}}$. The construction can be made in such a
way that edges do not intersect. Notice that condition {\bf 2.}~in the
definition of the set $T^M$ entails that a type-$2$ vertex is
connected by the preceding construction to exactly two type-$1$
vertices with the same label, so that we can erase all type-$2$
vertices. The resulting planar graph is a planar map. We view this map
as a rooted planar map by declaring that the distinguished edge is the
one corresponding to $k=0$, pointing from $\delta$, in the preceding
construction.

It follows from \cite{BdFG} that the preceding construction yields a
bijection $\Psi_{r}$ between $\ov{\Tr}_{M}$ and $\m_{r}$. Furthermore
it is not difficult to see that $\Psi_{r}$ satisfies the following two
properties~: let $(\tr,\e,\l)\in\ov{\Tr}_{M}$ and let
$\map=\Psi_{r}((\tr,\e,\l))$,
\begin{enumerate}
\item[(i)] the set $\f_\map$ is in one-to-one correspondence with the
  set $\tr^3\cup\tr^4$, more precisely, with every $v\in\tr^3$
  (resp.~$v\in\tr^4$) such that $\z_u(\tr)=(k,k',0,0)$ is associated a
  unique face of $\map$ whose degree is equal to $2k+k'+2$
  (resp.~$2k+k'+1$),
\item[(ii)] for every $l\geq1$, the set $\{a\in\v_\map:d(\dd,a)=l\}$
  is in one-to-one correspondence with the set $\{v\in\tr^1:\l_v=l\}$.
\end{enumerate}

\subsection{Boltzmann laws on multitype spatial trees}\label{secdeflois}

Let $\q$ be a regular critical weight sequence. We associate with $\q$
four probability measures on $\Z_+^4$ defined by~: \ba
\mu_\q^{(1)}(\{(0,0,k,0)\})&=&\frac{1}{Z_\q^+}\left(1-\frac{1}{Z_\q^+}\right)^k,\;\,k\geq0,\\
\mu_\q^{(2)}(\{(0,0,0,1)\})&=&1,\\
\mu_\q^{(3)}(\{(k,k',0,0)\})&=&\frac{(Z_\q^+)^k(Z_\q^\diamondsuit)^{k'}N_\bullet(k,k'){k+k'\choose k}q_{2+2k+k'}}{f_\q^\bullet(Z_\q^+,Z_\q^\diamondsuit)},\;\,k,k'\geq0,\\
\mu_\q^{(4)}(\{(k,k',0,0)\})&=&\frac{(Z_\q^+)^k(Z_\q^\diamondsuit)^{k'}N_\diamondsuit(k,k'){k+k'\choose
    k}q_{1+2k+k'}}{f_\q^\diamondsuit(Z_\q^+,Z_\q^\diamondsuit)},\;\,k,k'\geq0.
\ea We set
$\mug_\q=\left(\mu_\q^{(1)},\mu_\q^{(2)},\mu_\q^{(3)},\mu_\q^{(4)}\right)$
and ${\sf M}_{\mup_\q}=(m_{ij})_{1\leq i,j\leq4}$. The matrix ${\sf
  M}_{\mup_\q}$ is given by
$${\sf M}_{\mup_\q}=\left(\begin{array}{cccc} 0 & 0 & Z_\q^+-1&0\\
    0&0&0&1\\
    \frac{(Z_\q^+)^2}{Z_\q^+-1}\dd_xf_\q^\bullet(Z_\q^+,Z_\q^\diamondsuit)&\frac{Z_\q^+ Z_\q^{\diamondsuit}}{Z_\q^+-1}\dd_yf_\q^\bullet(Z_\q^+,Z_\q^\diamondsuit)&0&0\\
    \frac{Z_\q^+}{Z_\q^\diamondsuit}\,\dd_xf_\q^\diamondsuit(Z_\q^+,Z_\q^\diamondsuit)&\dd_yf_\q^\diamondsuit(Z_\q^+,Z_\q^\diamondsuit)&0&0 \\
\end{array}\right).$$
We see that ${\sf M}_{\mup_\q}$ is irreducible and has a spectral
radius $\varrho=1$. Thus $\mug_\q$ is critical. Let us denote by
$\ag=(a_1,a_2,a_3,a_4)$ the right eigenvector of ${\sf M}_{\mup_\q}$
with eigenvalue $1$ chosen so that $a_1+a_2+a_3+a_4=1$.

Let $\zetag_\q$ be the uniform ordering of $\mug_\q$. Note that if
$\w\in\W_4$ satisfies $w_j\in\{1,2\}$ for every
$j\in\{1,\ldots,|\w|\}$, then, by setting $k=p_1(\w)$ and
$k'=p_2(\w)$, we have \ba
\zeta_\q^{(3)}(\{\w\})&=&\frac{(Z_\q^+)^k(Z_\q^\diamondsuit)^{k'}N_\bullet(k,k')q_{2+2k+k'}}{f_\q^\bullet(Z_\q^+,Z_\q^\diamondsuit)},\\
\zeta_\q^{(4)}(\{\w\})&=&\frac{(Z_\q^+)^k(Z_\q^\diamondsuit)^{k'}N_\diamondsuit(k,k')q_{1+2k+k'}}{f_\q^\diamondsuit(Z_\q^+,Z_\q^\diamondsuit)}.
\ea

Let us now define a collection
$\nug=(\nu_{i,\w},i\in\{1,2,3,4\},\w\in\W_4)$ as follows.
\begin{enumerate}
\item[$\bullet$] For $i\in\{1,2\}$ the measure $\nu_{i,\w}$ is the
  Dirac mass at $\zero\in\R^{|\w|}$.
\item[$\bullet$] Let $\w\in\W_4$ be such that $p(\w)=(k,k',0,0)$. Then
  $\nu_{3,\w}$ is the distribution of the random vector
  $(X_1,X_1+X_2,\ldots,X_1+X_2+\ldots+X_{k+k'})$, where
  $(X_j+\ind{\{w_{j-1}=1\}},1\leq j\leq k+k'+1)$ (with $w_0=1$) is
  uniformly distributed on the set
$$A_{k,k'}=\left\{(n_1,\ldots,n_{k+k'})\in\Z_+^{k+k'+1}:n_1+\ldots+n_{k+k'+1}=k+1\right\}.$$
\item[$\bullet$] Let $\w\in\W_4$ be such that $p(\w)=(k,k',0,0)$. Then
  $\nu_{4,\w}$ is the distribution of the random vector
  $(X_1,X_1+X_2,\ldots,X_1+X_2+\ldots+X_{k+k'})$, where
  $(X_j+\ind{\{w_{j-1}=1\}},1\leq j\leq k+k'+1)$ (with $w_0=2$) is
  uniformly distributed on the set
$$B_{k,k'}=\left\{(n_1,\ldots,n_{k+k'})\in\Z_+^{k+k'+1}:n_1+\ldots+n_{k+k'+1}=k\right\}.$$
\item[$\bullet$] If $i\in\{3,4\}$ and if $\w\in\W_4$ does not satisfy
  $p_3(\w)=p_4(\w)=0$ then $\nu_{i,\w}$ is arbitrarily defined.
\end{enumerate}
Note that $\#A_{k,k'}=N_\bullet(k,k')$ and $\#B_{k,k'}=N_\diamondsuit(k,k')$.

Let us now introduce some notation. We have
$P^{(i)}_{\mup_\q}(\#\tr^1=n)>0$ for every $n\geq1$ and
$i\in\{1,2\}$. Then we may define, for every $n\geq1$, $i\in\{1,2\}$
and $x\in\R$, \ba
P_{\mup_\q}^{(i),n}(\d\tr\,\d\e)&=&P_{\mup_\q}^{(i)}\left(\d\tr\,\d\e\mid\#\tr^1=n\right),\\
\P^{(i),n}_{\mup_\q,\nup,x}(\d\tr\,\d\e\,\d\l)&=&\P^{(i)}_{\mup_\q,\nup,x}\left(\d\tr\,\d\e\,\d\l\mid\#\tr^1=n\right).
\ea Furthermore, we set for every $(\tr,\l,\e)\in\Tr^{(4)}$,
$$\un{\l}=\min\left\{\l_v:v\in\tr^1\setminus\{\varnothing\}\right\},$$ 
with the convention $\min\varnothing=\infty$. Finally we define for
every $n\geq1$, $i\in\{1,2\}$ and $x\geq0$, \ba
\ov{\P}^{(i)}_{\mup_\q,\nup,x}(\d\tr\,\d\e\,\d\l)&=&\P^{(i)}_{\mup_\q,\nup,x}(\d\tr\,\d\e\,\d\l\mid\un{\l}>0),\\
\ov{\P}^{(i),n}_{\mup_\q,\nup,x}(\d\tr\,\d\e\,\d\l)&=&\ov{\P}^{(i)}_{\mup_\q,\nup,x}\left(\d\tr\,\d\e\,\d\l\mid\#\tr^1=n\right).
\ea

The following proposition can be proved from Proposition 3 of
\cite{Mi2} in the same way as Corollary 2.3 of \cite{We}.

\begin{proposition}\label{loiimage}
  The probability measure $\B^r_\q(\cdot\mid\#\v_\map=n)$ is the image
  of $\,\ov{\P}^{(1),n}_{\mup_\q,\nup,1}$ under the mapping $\Psi_r$.
\end{proposition}

\section{A conditional limit theorem for multitype spatial trees}

Let $\q$ be a regular critical weight sequence. Recall from section
\ref{secdeflois} the definition of the offspring distribution
$\mug_\q$ associated with $\q$ and the definition of the spatial
displacement distributions $\nug$. To simplify notation we set
$\mug=\mug_\q$.

In view of applying a result of \cite{Mi1}, we have to take into
account the fact that the spatial displacements $\nug$ are not
centered distributions, and to this end we will need a shuffled
version of the spatial displacement distributions $\nug$. Let
$i\in[K]$ and $\w\in\W$. Set $n=|\w|$. We set
$\ovla{\w}=(w_n,\ldots,w_1)$ and we denote by $\ovla{\nu}_{i,\w}$ the
image of the measure $\nu_{i,\w}$ under the mapping
$S_n:(x_1,\ldots,x_n)\mapsto(x_n,\ldots,x_1)$. Last we set
$$\ovlra{\nu}_{i,\w}(\d\yg)=\frac{\nu_{i,\w}(\d\yg)+
  \ovla{\nu}_{i,\ovla{\w}}(\d\yg)}{2}.$$ We write
$\ovla{\nug}=(\ovla{\nu}_{i,\w},i\in[K],\w\in\W)$ and
$\ovlra{\nug}=(\ovlra{\nu}_{i,\w},i\in[K],\w\in\W)$.  

If $(\tr,\e,\l)$ is a multitype spatial tree, we denote by $C$ its
contour function and by $V$ its spatial contour function. Recall that
$C([0,1],\R)^2$ is equipped with the norm
$\|(f,g)\|=\|f\|_u\vee\|g\|_u$. The following result is a special case
of Theorem 4 in \cite{Mi1}.

\begin{theorem}\label{thconv} Let $\q$ be a regular critical weight
  sequence. There exists two scaling constants $\Ad_\q>0$ and
  $\Bd_\q>0$ such that for $i\in\{1,2\}$, the law under
  $\P_{\mup,\ovlra\nup,0}^{(i),n}$ of
$$\left(\left(\Ad_\q\,\frac{C(2(\#\tr-1)s)}{n^{1/2}}\right)_{0\leq s\leq1},\left(\Bd_\q\,\frac{V(2(\#\tr-1)s)}{n^{1/4}}\right)_{0\leq s\leq1}\right)$$
converges as $n\to\infty$ to the law of $(\b,\r^0)$. The convergence
holds in the sense of weak convergence of probability measures on
$C([0,1],\R)^2$.
\end{theorem}
Note that Theorem 4 in \cite{Mi1} deals with the so-called
height process instead of the contour process. However, we can deduce
Theorem \ref{thconv} from \cite{Mi1} by classical arguments (see
e.g.~\cite{LGDEA}). Moreover, the careful reader will notice that the
spatial displacements $\ovlra{\nug}$ depicted above are not all
centered, and thus may compromise the application of \cite[Theorem
4]{Mi1}. However, it is explained in \cite[Sect.\ 3.3]{Mi2} how a simple
modification of these laws can turn them into centered distributions,
by appropriate translations. More precisely, one can couple the
spatial trees associated with $\ovlra{\nug}$ and its centered version
so that the labels of vertices differ by at most $1/2$ in absolute
value, which of course does not change the limiting behavior of the
label function rescaled by $n^{-1/4}$.

In this section, we will prove a conditional version of Theorem
\ref{thconv}. Before stating this result, we establish a corollary of
Theorem \ref{thconv}. To this end we set \ba
Q_{\mup}(\d\tr\,\d\e)&=&P_{\mup}^{(1)}(\d\tr\,\d\e\mid c_\varnothing(\tr)=1),\\
\Q_{\mup,\ovlra\nup}(\d\tr\,\d\e\,\d\l)&=&\P_{\mup,\ovlra\nup,0}^{(1)}
(\d\tr\,\d\e\,\d\l\mid c_\varnothing(\tr)=1).  \ea Notice that this
conditioning makes sense since $\mu^{(1)}(\{(0,0,1,0)\})>0$. We may
also define for every $n\geq1$, \ba
Q_{\mup}^n(\d\tr\,\d\e)&=&Q_{\mup}\left(\d\tr\,\d\e\mid\#\tr^1=n\right),\\
\Q_{\mup,\ovlra\nup}^n(\d\tr\,\d\e\,\d\l)&=&\Q_{\mup,\ovlra\nup}\left(\d\tr\,\d\e\,\d\l\mid\#\tr^1=n\right).
\ea The following corollary can be proved from Theorem \ref{thconv} in
the same way as Corollary 2.2 in \cite{We}.

\begin{corollary}\label{corconvQ}
  Let $\q$ be a regular critical weight sequence. The law under
  $\Q_{\mup,\ovlra\nup}^n$ of
$$\left(\left(\Ad_\q\,\frac{C(2(\#\tr-1)s)}{n^{1/2}}\right)_{0\leq s\leq1},\left(\Bd_\q\,\frac{V(2(\#\tr-1)s)}{n^{1/4}}\right)_{0\leq s\leq1}\right)$$
converges as $n\to\infty$ to the law of $(\b,\r^0)$. The convergence
holds in the sense of weak convergence of probability measures on
$C([0,1],\R)^2$.
\end{corollary}

Recall from section \ref{serpent} that $(\ov{\b}^0,\ov{\r}^0)$ denotes
the conditioned Brownian snake with initial point $0$.

\begin{theorem}\label{thconvcond} Let $\q$ be a regular critical
  weight sequence. For every $x\geq0$, the law under
  $\ov{\P}_{\mu^\q,\ovlra\nu,x}^{(1),n}$ of
$$\left(\left(\Ad_\q\,\frac{C(2(\#\tr-1)s)}{n^{1/2}}\right)_{0\leq s\leq1},\left(\Bd_\q\,\frac{V(2(\#\tr-1)s)}{n^{1/4}}\right)_{0\leq s\leq1}\right)$$
converges as $n\to\infty$ to the law of $(\ov{\b}^0,\ov{\r}^0)$. The
convergence holds in the sense of weak convergence of probability
measures on $C([0,1],\R)^2$.
\end{theorem}

In the same way as in the proof of Theorem 3.3 in \cite{We}, we will
follow the lines of the proof of Theorem 2.2 in \cite{LG} to prove
Theorem \ref{thconvcond}.

\subsection{Rerooting spatial trees}\label{sec:rero-spat-trees}

If $(\tr,\e)\in T_M$, we write $\dd\tr=\{u\in \tr: c_u(\tr)=0\}$ for
the set of all leaves of $\tr$, and we write
$\dd_1\tr=\dd\tr\cap\tr^1$ for the set of leaves of $\tr$ which are of
type $1$. Let $w_0\in\tr$. Recall that
$\u^\ast=\u\setminus\{\varnothing\}$. We set
$$\tr^{(w_0)}=\tr\setminus\left\{w_0u\in\tr:u\in\u^\ast\right\},$$ 
and we write $\e^{(w_0)}$ for the restriction of the funtion $\e$ to
the truncated tree $\tr^{(w_0)}$.

Let $v_0=u^1\ldots u^{2p}\in\u^\ast$ and $(\tr,\e)\in T_M$ such that
$v_0\in\tr^1$. We define $k=k(v_0,\tr)$ and $l=l(v_0,\tr)$ in the
following way. Write $\xi=\#\tr-1$ and $u_0,u_1,\ldots,u_{2\xi}$ for
the search-depth sequence of $\tr$. Then we set \ba
k&=&\min\{i\in\{0,1,\ldots,2\xi\}:u_i=v_0\},\\
l&=&\max\{i\in\{0,1,\ldots,2\xi\}:u_i=v_0\}, \ea which means that $k$
is the time of the first visit of $v_0$ in the evolution of the
contour of $\tr$ and that $l$ is the time of the last visit of
$v_0$. Note that $l\geq k$ and that $l=k$ if and only if
$v_0\in\dd\tr$. For every $s\in[0,2\xi-(l-k)]$, we set
$$\wh{C}^{(v_0)}(s)=C(k)+C(\llbr k-s\rrbr)-2\inf_{u\in[k\wedge\llbr k-s\rrbr,k\vee\llbr k-s\rrbr]}C(u),$$
where $C$ is the contour function of $\tr$ and $\llbr k-s\rrbr$ stands
for the unique element of $[0,2\xi)$ such that $\llbr
k-s\rrbr-(k-s)=0$ or $2\xi$. Then there exists a unique plane tree
$\wh{\tr}^{(v_0)}$ whose contour function is
$\wh{C}^{(v_0)}$. Informally, $\wh{\tr}^{(v_0)}$ is obtained from
$\tr$ by removing all vertices that are descendants of $v_0$, by
re-rooting the resulting tree at $v_0$, and finally by reversing the
planar orientation. Furthermore we see that $\wh{v}_0=1u^{2p}\ldots
u^2$ belongs to $\wh{\tr}^{(v_0)}$. In fact, $\wh{v}_0$ is the vertex
of $\wh{\tr}^{(v_0)}$ corresponding to the root of the initial
tree. At last notice that $c_\root(\,\wh{\tr}^{(v_0)})=1$.

We now define the function $\wh{\e}^{(v_0)}$. To this end, for $u\in
\llbr\root, v_0\rrbr\setminus\{v_0\}$, let
$j(u,v_0)\in\{1,\ldots,c_u(\tr)\}$ be such that $uj(u,v_0)\in
\llbr\root,v_0\rrbr$. Then set
\ba
\llbr\root,v_0\rrbr_2^3&=&\left\{u\in\llbr\root,v_0\rrbr\cap\tr^3:
  \e(uj(u,v_0))=2\right\}\\
\llbr\root,v_0\rrbr_1^4&=&\left\{u\in\llbr\root,v_0\rrbr\cap\tr^4:
  \e(uj(u,v_0))=1\right\}.  \ea 
For every $u\in\wh{\tr}^{(v_0)}$, we denote by $\ov{u}$ the vertex which
corresponds to $u$ in the tree $\tr$. We then set $\wh{\e}^{(v_0)}(u)=\e(\ov{u})$, except in the following cases~: 
\begin{equation}
\left\{\begin{array}{l}
\mbox{if }\ov{u}\in\llbr\root,v_0\rrbr_3^2 
\mbox{ then }\wh{\e}^{(v_0)}(u)=4, \\
\mbox{if }\ov{u}\in\llbr\root,v_0\rrbr_1^4 
\mbox{ then }\wh{\e}^{(v_0)}(u)=3.
\label{discussion}
\end{array}\right. 
\end{equation}
Since $v_0\in\tr^1$ we have
$\#\llbr\root,v_0\rrbr_2^3=\#\llbr\root,v_0\rrbr_1^4$. Indeed, if
$1=e_0,e_1,\ldots,e_{2p}=1$ is the sequence of types of elements of
$\llbr\root,v_0\rrbr$ listed according to their generations, then this
list is a concatenation of patterns of the form $13\ov{24}1$, where by
$\ov{24}$ we mean an arbitrary (possibly empty) repetition of the
pattern $24$. If at least one $24$ occurs then the second and
antepenultimate element of the pattern $13\ov{24}1$ correspond
respectively to exactly one element of $\llbr\root,v_0\rrbr_2^3$ and
$\llbr\root,v_0\rrbr_1^4$, while no term of a pattern $131$
corresponds to such elements.

Notice that if $(\,\wh{\tr}^{(v_0)},\wh{\e}^{(v_0)})=(\t,\ed)$, then
$(\tr^{(v_0)},\e^{(v_0)})=(\,\wh{\t}^{(\wh{v}_0)},\wh{\ed}^{(\wh{v}_0)})$.
Moreover, if $u\in\t\setminus\{\root,\wh{v}_0\}$ then we have
$c_u(\t)=c_{\ov{u}}(\wh{\t}^{(\wh{v}_0}))$. Recall that if $\w=(w_1,\ldots,w_n)$ we write $\ovla{\w}=(w_n,\ldots,w_1)$. To be more accurate, it holds that $\w_u(\t)=\overleftarrow{\w}_{\ov{u}}(\wh{\t}^{(\wh{v}_0)})$ except in the following cases~: 
\be\left\{\begin{array}{l}
\mbox{if }\ov{u}\in\llbr\root,v_0\rrbr\setminus(\llbr\root,v_0\rrbr_2^3\cap\llbr\root,v_0\rrbr_1^4) \mbox{ then }\w_u(\t)=\overleftarrow{\w}^{j(\ov{u},v_0),\e(\ov{u})}_{\ov{u}}(
  \wh{\t}^{(\wh{v}_0)}), \\
\mbox{if }\ov{u}\in\llbr\root,v_0\rrbr_2^3 \mbox{ then }
  \w_u(\t)=\overleftarrow{\w}^{j(\ov{u},v_0),1}_{\ov{u}}(
  \wh{\t}^{(\wh{v}_0)}), \\
\mbox{if }\ov{u}\in\llbr\root,v_0\rrbr_1^4 \mbox{ then }
  \w_u(\t)=\overleftarrow{\w}^{j(\ov{u},v_0),2}_{\ov{u}}(
\wh{\t}^{(\wh{v}_0)}),
\end{array}\right.\label{discussion2}
\ee
where for $\w\in\W$, $n=|\w|$, and $1\leq j\leq n$, we set
$$
\left\{\begin{array}{l}
  \w^{j,1}=(w_{j+1},\ldots,w_{n},1,w_{1},\ldots,w_{j-1}),\\
  \w^{j,2}=(w_{j+1},\ldots,w_{n},2,w_{1},\ldots,w_{j-1}).
\end{array}\right.
$$
In particular, if $\ov{u}\in\llbr\root,v_0\rrbr_2^3$ (resp.\ $\llbr\root,v_0\rrbr_1^4$) with $p(\w_{\ov{u}}(\wh{\t}^{(\wh{v}_0)}))=(k,k',0,0)$ then \linebreak $p(\w_u(\t))=(k+1,k'-1,0,0)$ (resp.\ $(k-1,k'+1,0,0)$), while $p(\w_{\ov{u}}(\wh{\t}^{(\wh{v}_0)}))=p(\w_u(\t))$ otherwise.

\begin{figure}[htbp]
  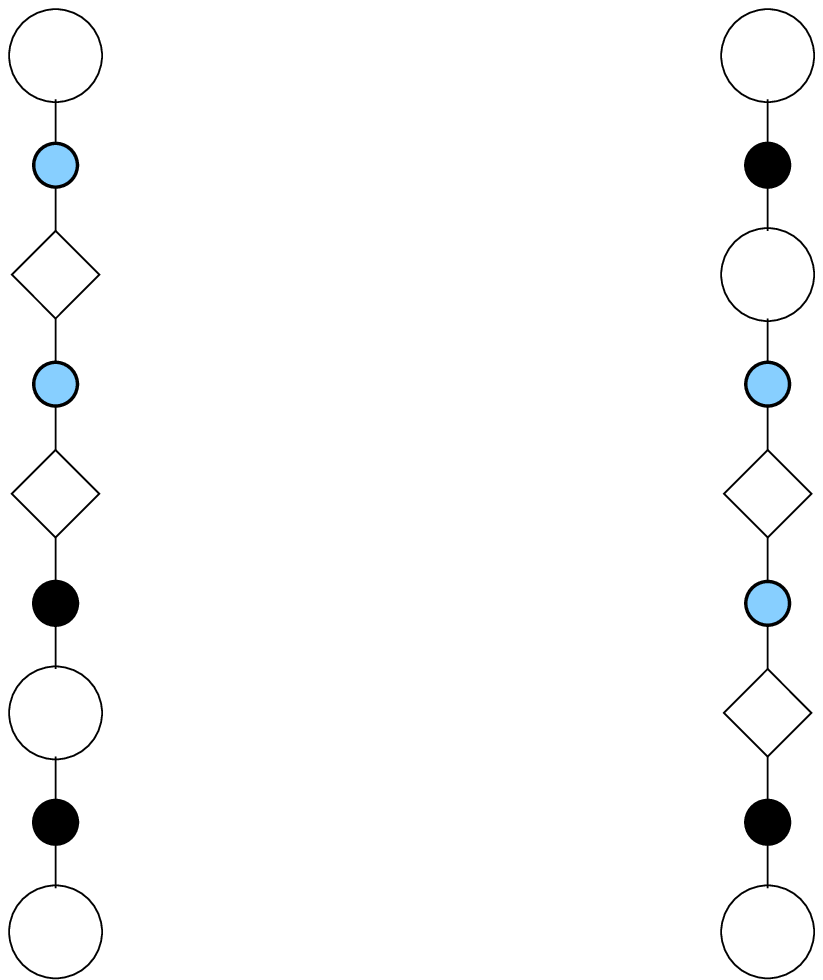
\caption{The branch leading from $\varnothing$ to $v_0$, and the
  corresponding branch in the tree $\widehat{\bt}^{(v_0)}$: the branch
  is put upside-down and the vertices of $\llbr\root,v_0\rrbr_2^3$ and
  $\llbr\root,v_0\rrbr_1^4$ interchange their roles.}
\end{figure}

Recall the definition of the probability measure $Q_\mup$.

\begin{lemma}\label{reenr}
  Let $v_0\in\u^\ast$ be of the form $v_0=1u^2\ldots u^{2p}$ for some
  $p\in\N$. Assume that
$$Q_\mup\left(v_0\in\tr^1\right)>0.$$ 
Then the law of the re-rooted multitype tree
$(\,\wh{\tr}^{(v_0)},\wh{\e}^{(v_0)})$ under $Q_\mup(\cdot\mid
v_0\in\tr^1)$ coincides with the law of the multitype tree
$(\tr^{(\wh{v}_0)},\e^{(\wh{v}_0)})$ under
$Q_\mup(\cdot\mid\wh{v}_0\in\tr^1)$.
\end{lemma}

\proof Let $(\t,\ed)\in T_M$ such that $\wh{v}_0\in\dd_1\t$. We have
$$Q_\mup\left((\,\wh{\tr}^{(v_0)},\wh{\e}^{(v_0)})=(\t,\ed)\right)=
Q_\mup\left((\tr^{(v_0)},\e^{(v_0)})=(\wh{\t}^{(\wh{v}_0)},
  \wh{\ed}^{(\wh{v}_0)})\right).$$ And
$$Q_\mup\left((\tr^{(v_0)},\e^{(v_0)})=(\wh{\t}^{(\wh{v}_0)},
  \wh{\ed}^{(\wh{v}_0)})\right)=\prod_{\ov{u}\in\wh{\t}^{(\wh{v}_0)}
  \setminus\{\root,v_0\}}\zeta^{(\wh{\ed}^{(\wh{v}_0)}(u))}(
\w_u(\wh{\t}^{(\wh{v}_0)})),$$
$$Q_\mup\left((\tr^{(\wh{v}_0)},\e^{(\wh{v}_0)})=(\t,\ed)\right)=
\prod_{u\in\t\setminus\{\root,\wh{v}_0\}}\zeta^{(\ed(u))}(\w_u(\t)).$$
By the above discussion around (\ref{discussion2}), the terms corresponding to $u,\ov{u}$ in
these two products are all equal, except for those corresponding to
vertices
$\ov{u}\in\llbr\root,v_0\rrbr_3^2\cup\llbr\root,v_0\rrbr_4^1$.

Let $k\geq0$ and $k'\geq1$. We have
$N_\diamondsuit(k+1,k'-1)=N_\bullet(k,k')$ which implies that \ba
\frac{\mu^{(4)}(k+1,k'-1,0,0)}{{k+k'\choose k+1}}&=&
\frac{(Z_\q^+)^{k+1}(Z_\q^\diamondsuit)^{k'-1}
  N_\diamondsuit(k+1,k'-1)q_{1+2(k+1)+k'-1}}{f_\q^\diamondsuit(
  Z_\q^+,Z_\q^\diamondsuit)}\\
&=&\frac{Z_\q^+f_\q^\bullet(Z_\q^+,Z_\q^\diamondsuit)}{
  Z_\q^\diamondsuit f_\q^\diamondsuit(Z_\q^+,Z_\q^\diamondsuit)}\,
\frac{\mu^{(3)}(k,k',0,0)}{{k+k'\choose k}}\\
&=&\frac{Z_\q^+-1}{(Z_\q^\diamondsuit)^2}\,\frac{\mu^{(3)}(k,k',0,0)}{
  {k+k'\choose k}}.  \ea Likewise let $k\geq1$ and $k'\geq0$. We have
$N_\bullet(k-1,k'+1)=N_\diamondsuit(k,k')$ which implies that \ba
\frac{\mu^{(3)}(k-1,k'+1,0,0)}{{k+k'\choose k-1}}&=&\frac{(
  Z_\q^+)^{k-1}(Z_\q^\diamondsuit)^{k'+1}N_\bullet(k-1,k'+1)
  q_{2+2(k-1)+k'+1}}{f_\q^\bullet(Z_\q^+,Z_\q^\diamondsuit)}\\
&=&\frac{Z_\q^\diamondsuit
  f_\q^\diamondsuit(Z_\q^+,Z_\q^\diamondsuit)}{
  Z_\q^+f_\q^\bullet(Z_\q^+,Z_\q^\diamondsuit)}\,
\frac{\mu^{(4)}(k,k',0,0)}{{k+k'\choose k}}\\
&=&\frac{(Z_\q^\diamondsuit)^2}{Z_\q^+-1}\,\frac{\mu^{(4)}(k,k',0,0)}{
  {k+k'\choose k}}.  \ea Using the relation between
$p(\w_{\ov{u}}(\wh{\t}^{(\wh{v}_0)}))$ and $p(\w_u(\t))$ discussed above
for elements of $\llbr\root,v_0\rrbr_2^3\cup\llbr\root,v_0\rrbr_1^4$,
we obtain \ba &&\hspace{-1.5cm}Q_\mup\left((\tr^{(v_0)},\e^{(v_0)})=
  (\wh{\t}^{(\wh{v}_0)},\wh{\ed}^{(\wh{v}_0)})\right)\\
&=&\left(\frac{Z_\q^+-1}{(Z_\q^\diamondsuit)^2}\right)^{
  \#\llbr\root,v_0\rrbr_2^3-\#\llbr\root,v_0\rrbr_1^4}Q_\mup
\left((\tr^{(\wh{v}_0)},\e^{(\wh{v}_0)})=(\t,\ed)\right)\\
&=&Q_\mup\left((\tr^{(\wh{v}_0)},\e^{(\wh{v}_0)})=(\t,\ed)\right), \ea
implying that
\be\label{eqreroot}
Q_\mup\left((\,\wh{\tr}^{(v_0)},\wh{\e}^{(v_0)})=(\t,\ed)\right)=
Q_\mup\left((\tr^{(\wh{v}_0)},\e^{(\wh{v}_0)})=(\t,\ed)\right).\ee
To conclude the proof we use (\ref{eqreroot}) to get that
\ba
Q_\mup(v_0\in\tr)&=&\sum_{\{(\t,\ed)\in T_M:v_0\in\dd_1\t\}}\Q_\mup\left((\tr^{(v_0)},\e^{(v_0)})=(\t,\ed)\right)\\
&=&\sum_{\{(\t,\ed)\in T_M:v_0\in\dd_1\t\}}\Q_\mup\left((\wh{\tr}^{(v_0)},\wh{\e}^{(v_0)})=(\wh{\t}^{(v_0)},\wh{\ed}^{(v_0)})\right)\\
&=&\sum_{\{(\t,\ed)\in T_M:v_0\in\dd_1\t\}}\Q_\mup\left((\tr^{(\wh{v}_0)},\e^{(\wh{v}_0)})=(\wh{\t}^{(v_0)},\wh{\ed}^{(v_0)})\right)\\
&=&\sum_{\{(\t',\ed')\in T_M:\wh{v}_0\in\dd_1\t'\}}\Q_\mup\left((\tr^{(\wh{v}_0)},\e^{(\wh{v}_0)})=(\t',\ed')\right)\\
&=&Q_\mup(\wh{v}_0\in\tr).
\ea
\cq

If $(\tr,\e,\l)\in\Tr_M$ and $v_0\in\tr^1$, the re-rooted multitype
spatial tree $(\wh{\tr}^{(v_0)},\wh{\e}^{(v_0)},\wh{\l}^{(v_0)})$ is
defined as follows. If $u\in\wh{\tr}^{(v_0)}$, recall that $\ov{u}$
denotes the vertex which correponds to $u$ in the tree $\tr$ and that
$\check{u}$ denotes its father (in the tree $\wh{\tr}^{(v_0)}$).
\begin{enumerate}
\item[$\bullet$] If $\wh{\e}^{(v_0)}(u)\in\{1,2\}$ then
  $\wh{\l}^{(v_0)}_u=\l_{\ov{u}}-\l_{v_0}$.
\item[$\bullet$] If $\wh{\e}^{(v_0)}(u)\in\{3,4\}$ then
  $\wh{\l}^{(v_0)}_u=\l^{(v_0)}_{\check{u}}$.
\end{enumerate}
Let $n=c_u(\wh{\tr}^{(v_0)})$. Observe that when $\ov{u}\notin\llbr \varnothing,v_0\rrbr$, then the spatial displacements between $u$ and its offspring is left unchanged by the re-rooting, meaning that
$$\left(\wh{\l}^{(v_0)}_{ui}-\wh{\l}^{(v_0)}_u,1\leq i\leq n\right)=\left(\l^{(\wh{v}_0)}_{\ov{u}i}-\l^{(\wh{v}_0)}_{\ov{u}},1\leq i\leq n\right)\, .$$
Otherwise, if $\ov{u}\in\llbr\varnothing,v_0\rrbr$, set $j=j(u,v_0)$, and define the mapping
  $$\phi_{n,j}:(x_1,\ldots,x_n)\mapsto(x_{j-1}-x_j,\ldots,x_1-x_j,-x_j,
  x_n-x_j,\ldots,x_{j+1}-x_j)\, .$$ Then observe that the spatial displacements are affected in the following way:
\begin{equation}\label{discussion4}
\left(\wh{\l}^{(v_0)}_{ui}-\wh{\l}^{(v_0)}_u,1\leq i\leq n\right)=\phi_{n,j}\left(\left(\l^{(\wh{v}_0)}_{\ov{u}i}-\l^{(\wh{v}_0)}_{\ov{u}},1\leq i\leq n\right)\right)\, .
\end{equation}

If $(\tr,\e,\l)\in\Tr_M$ and $w_0\in\tr$, we also consider the
multitype spatial tree $(\tr^{(w_0)},\e^{(w_0)},\l^{(w_0)})$ where
$\l^{(w_0)}$ is the restriction of $\l$ to the tree $\tr^{(w_0)}$.

Recall the definition of the probability measure
$\Q_{\mup,\ovlra{\nup}}$.

\begin{lemma}\label{reenrspat}
  Let $v_0\in\u^\ast$ be of the form $v_0=1u^2\ldots u^{2p}$ for some
  $p\in\N$. Assume that
$$Q_\mup\left(v_0\in\tr^1\right)>0.$$ 
Then the law of the re-rooted multitype spatial tree
$(\,\wh{\tr}^{(v_0)},\wh{\e}^{(v_0)},\wh{\l}^{(v_0)})$ under the
measure $\Q_{\mup,\ovlra{\nup}}(\cdot\mid v_0\in\tr^1)$ coincides with
the law of the multitype spatial tree
$(\tr^{(\wh{v}_0)},\e^{(\wh{v}_0)},\l^{(\wh{v}_0)})$ under the measure
$\Q_{\mup,\ovlra{\nup}}(\cdot\mid\wh{v}_0\in\tr^1)$.
\end{lemma}

This lemma is a simple consequence of Lemma \ref{reenr} and our observations around (\ref{discussion4}) on the spacial displacements 
$\ovlra{\nug}$,
combined with the discussion of how the set of children of various
vertices are affected by re-rooting, see (\ref{discussion}) and (\ref{discussion2}).

\begin{lemma}
  Let $\w\in\W$ such that $p_3(\w)=p_4(\w)=0$. Set $n=|\w|$ and let
  $j\in\{1,\ldots,n\}$. 
  \begin{enumerate}
\item[(i)] The image of the measure $\ovlra{\nu}_{3,\w}$ under the
  mapping $\phi_{n,j}$ is
  \begin{enumerate}
\item[(a)] the measure $\ovlra{\nu}_{3,\w^{j,1}}$ if $w_j=1$,
\item[(b)] the measure $\ovlra{\nu}_{4,\w^{j,2}}$ if $w_j=2$.
\end{enumerate}
\item[(ii)] The image of the measure $\ovlra{\nu}_{4,\w}$ under the
  mapping $\phi_{n,j}$ is
  \begin{enumerate}
\item[(a)] the measure $\ovlra{\nu}_{3,\w^{j,1}}$ if $w_j=1$,
\item[(b)] the measure $\ovlra{\nu}_{4,\w^{j,2}}$ if $w_j=2$.
\end{enumerate}
\end{enumerate}
\end{lemma}

\proof We first suppose that $w_j=1$. Set $k=p_1(\w)$, $k'=p_2(\w)$
and $w_0=1$. Define $\wt{\phi}_{n,j}=S_n\circ\phi_{n,j}$, where as before 
$S_n$ stands for the mapping $(x_1,\ldots,x_n)\mapsto(x_n,\ldots,x_1)$. 
We consider a uniform vector
$(X_{l}+\ind{\{w_{l-1}=1\}},1\leq l\leq n+1)$ on the set $A_{k,k'}$
and we set $X^{(j)}=(X_{j+1},\ldots,X_{n+1},X_1,\ldots,X_j)$. Since
$w_0=w_j=1$, the vector
$$\left(X^{(j)}_{l}+\ind{\{w^{j,1}_{l-1}=1\}},1\leq l\leq n+1\right)$$ 
is uniformly distributed on the set $A_{k,k'}$ and the measure
$\nu_{3,\w^{j,1}}$ is the law of the vector
$$\left(X^{(j)}_1,X^{(j)}_1+X^{(j)}_2,\ldots,X^{(j)}_1+\ldots+X^{(j)}_n\right).$$
Furthermore we notice that $X_1+X_2+\ldots+X_{n+1}=0$. This implies
that
$$\left(X^{(j)}_1,X^{(j)}_1+X^{(j)}_2,\ldots,X^{(j)}_1+\ldots+X^{(j)}_n\right)=\wt{\phi}_{n,j}(X_1,X_1+X_2,\ldots,X_1+\ldots+X_n),$$
which means that the measure $\nu_{3,\w^{j,1}}$ is the image of
$\nu_{3,\w}$ under the mapping $\wt{\phi}_{n,j}$. Since
$\wt{\phi}_{n,j}\circ S_n=S_n\circ\phi_{n,n-j+1}$, we obtain together with
what precedes that the measure $\ovla{\nu}_{3,\ovla{\w}^{n-j+1,1}}$ is
the image of $\ovla{\nu}_{3,\ovla{\w}}$ under the mapping
$\wt{\phi}_{n,j}$. Thus $\ovlra{\nu}_{3,\w^{j,1}}$ is the image of $\ovlra{\nu}_{3,\w}$ under the mapping $\wt{\phi}_{n,j}$. Hence, it is the image of the same measure under $\phi_{n,j}$, being invariant under the action of $S_n$. Thus we get the first part of the lemma. The other
assertions can be proved in the same way. \cq

If $(\tr,\e,\l)\in\Tr_M$, we set
$$\un{\l}=\min\left\{\l_v:v\in\tr^1\setminus\{\root\}\right\},$$
$$\Delta_1=\left\{v\in\tr^1:\l_v=\min\left\{\l_w:w\in\tr^1\right\}\right\}.$$
We also denote by $v_m$ the first element of $\Delta_1$ in the
lexicographical order.

The following two Lemmas can be proved from Lemma \ref{reenrspat} in
the same way as Lemma 3.3 and Lemma 3.4 in \cite{LG}.

\begin{lemma}\label{reenrmin1}
  For any nonnegative measurable functional $F$ on $\Tr_M$,
$$\Q_{\mup,\ovlra{\nup}}\left(F\left(\wh{\tr}^{(v_m)},\wh{\e}^{(v_m)},\wh{\l}^{(v_m)}\right)\ind{\{\#\Delta_1=1,v_m\in\dd_1\tr\}}\right)=\Q_{\mup,\ovlra{\nup}}\left(F(\tr,\e,\l)(\#\dd_1\tr)\ind{\{\un{\lp}>0\}}\right).$$
\end{lemma}

\begin{lemma}\label{reenrmin2}
For any nonnegative measurable functional $F$ on $\Tr_M$,
$$\Q_{\mup,\ovlra{\nup}}\left(\sum_{v_0\in\Delta_1\cap\dd_1\tr}F\left(\wh{\tr}^{(v_0)},\wh{\e}^{(v_0)},\wh{\l}^{(v_0)}\right)\right)=\Q_{\mup,\ovlra{\nup}}\left(F(\tr,\e,\l)(\#\dd_1\tr)\ind{\{\un{\lp}\geq0\}}\right).$$
\end{lemma}

\subsection{Estimates for the probability of staying on the positive side}\label{sec:estpos}

In this section we will derive upper and lower bounds for the
probability $\P^n_{\mup,\ovlra\nup,x}(\un{\l}>0)$ as $n\to\infty$. We
first state a lemma which is a direct consequence of Lemma 6 in
\cite{Mi1}.

\begin{lemma}\label{asympt} There exist two constants $c_0>0$ and
  $c_1>0$ such that
  \ba
  n^{3/2}P_\mup\left(\#\tr^1=n\right)&\build{\la}_{n\to\infty}^{}&c_0,\\
  n^{3/2}Q_\mup\left(\#\tr^1=n\right)&\build{\la}_{n\to\infty}^{}&c_1.
  \ea
\end{lemma}

We now establish a preliminary estimate concerning the number of
leaves of type $1$ in a tree with $n$ vertices of type $1$. Write $\zero\in\R^4$.

\begin{lemma}\label{feuilles}
  There exists a constant $\beta>0$ such that for all $n$ sufficiently
  large,
$$P_\mup\left(|\#\dd_1\tr-\mu^{(1)}(\{\zero\})n|>n^{3/4},\#\tr^1=n\right)\leq e^{-n^{\beta}}.$$
\end{lemma}

\proof Let $(\tr,\e)\in T_M$. Recall that $\xi=\#\tr-1$. Let $v(0)=\root\prec v(1)\prec\ldots\prec v(\xi)$ be the vertices of $\tr$ listed in
lexicographical order. For every $n\in\{0,1,\ldots,\xi\}$, we define
$R_n=(R_n(k))_{k\geq1}$ as follows. For every
$k\in\{1,\ldots,|v(n)|\}$, we write $v(n,k)$ for the ancestor of
$v(n)$ at generation $k$ and we let
$$v_1(n,k)\prec\ldots\prec v_m(n,k)$$
be the younger brothers of $v(n,k)$ listed in lexicographical
order. Here younger brothers are those brothers which have not yet
been visited at time $n$ in search-depth sequence. Then we set
$$R_n(k)=\left(\e(v_1(n,k)),\ldots,\e(v_m(n,k))\right)$$
if $m\geq1$ and $R_n(k)=\varnothing$ if $m=0$. For every $k>|v(n)|$,
we set $R_n(k)=\varnothing$. By abuse of notations we assimilate $R_n$
with $(R_n(1),\ldots,R_n(|v(n)|))$ and let $R_n=\varnothing$ if
$|v(n)|=0$. Standard arguments (see e.g.~\cite{LGLJ} for similar
results) show that $(R_n,\e(v(n)),|v(n)|)_{0\leq n\leq\xi}$ has the
same distribution as $(R'_n,e'_n,h'_n)_{0\leq n\leq T'-1}$, where
$(R'_n,e'_n,h'_n)_{n\geq0}$ is a Markov chain starting at
$(\varnothing,1,0)$, whose transition probabilities are specified as
follows~:

\begin{enumerate}
\item[$\bullet$] $((\r_1,\ldots,\r_h),i,h) \to
  ((\r_1,\ldots,\r_h,\r_{h+1}^+),\r_{h+1}(1),h+1)$ with probability
  $\zetag^{(i)}(\{\r_{h+1}\})$ where
  $\r^+_{h+1}=(\r_{h+1}(2),\ldots,\r_{h+1}(|\r_{h+1}|))$, for
  $\r_{h+1}\in\W_4\setminus\{\varnothing\}$, $i\in\{1,2,3,4\}$,
  $\r_1,\ldots,\r_h\in\W_4$ and $h\geq0$,
\item[$\bullet$] $((\r_1,\ldots,\r_h),i,h)\to
  ((\r_1,\ldots,\r_{k-1},\r^+_k),\r_k(1),k)$ with probability
  $\zetag^{(i)}(\{\varnothing\}),$ whenever $h\geq 1$ and
  $\{m\geq1:\r_m\neq\varnothing\}\neq\varnothing$, and where
  $k=\sup\{m\geq1:\r_m\neq\varnothing\}$, for $i\in\{1,2,3,4\}$,
  $\r_1,\ldots,\r_h\in\W_4$,  
\item[$\bullet$] $((\varnothing,\ldots,\varnothing),i,h)\to
  (\varnothing,1,0)$ with probability $\zetag^{(i)}(\{\varnothing\})$
  for $i\in\{1,2,3,4\}$ and $h\geq0$,
\end{enumerate}
and finally
$$T'=\inf\left\{n\geq1:
h'_n
=0
\right\}.$$
Write $\PP'$ for the probability measure under which
$(R'_n,e'_n,h'_n)_{n\geq0}$ is defined. We define a sequence of
stopping times $(\tau'_j)_{j\geq0}$ by $\tau'_0=0$ and
$\tau'_{j+1}=\inf\{n>\tau'_j:e'_n=1\}$ for every $j\geq0$. At last we
set for every $j\geq0$,
$$X'_j=\ind{}\left\{h_{\tau'_j+1}\leq h_{\tau'_j}\right\}.$$
Thus we have, \ba &&\hspace{-1cm}P_\mup\left(|\#\dd_1\tr-
  \mug^{(1)}(\{\zero\})n|>n^{3/4},\#\tr^1=n\right)\\
&=&\PP'\left(\Big|\sum_{j=0}^{n-1}X'_j-
  \mug^{(1)}(\{\zero\})n\Big|>n^{3/4},\tau'_{n-1}<
  T'\leq\tau'_n
\right)\\
&\leq&\PP'\left(\Big|\sum_{j=0}^{n-1}X'_j-\mug^{(1)}(\{\zero\})n
  \Big|>n^{3/4}
\right).  \ea Thanks to the strong Markov property, under the
probability measure $\PP'(\cdot\mid e'_0=1)$, the random variables
$X'_j$ are independent and distributed according to the Bernoulli
distribution with parameter
$\zetag^{(1)}(\{\varnothing\})=\mug^{(1)}(\{\zero\})$. So we get the
result using a standard moderate deviations inequality and Lemma
\ref{asympt}.\cq

We will now state a lemma which plays a crucial role in the proof of
the main result of this section. To this end, recall the definition of
$v_m$ and the definition of the probability measure
$\Q^n_{\mup,\ovlra{\nup}}$.

\begin{lemma}\label{argmin}
There exists a constant $c>0$ such that for all $n$ sufficiently large,
$$\Q_{\mup,\ovlra{\nup}}^n(v_m\in\dd_1\tr)\geq c.$$
\end{lemma}

\begin{figure}[htbp]
\includegraphics[scale=.7]{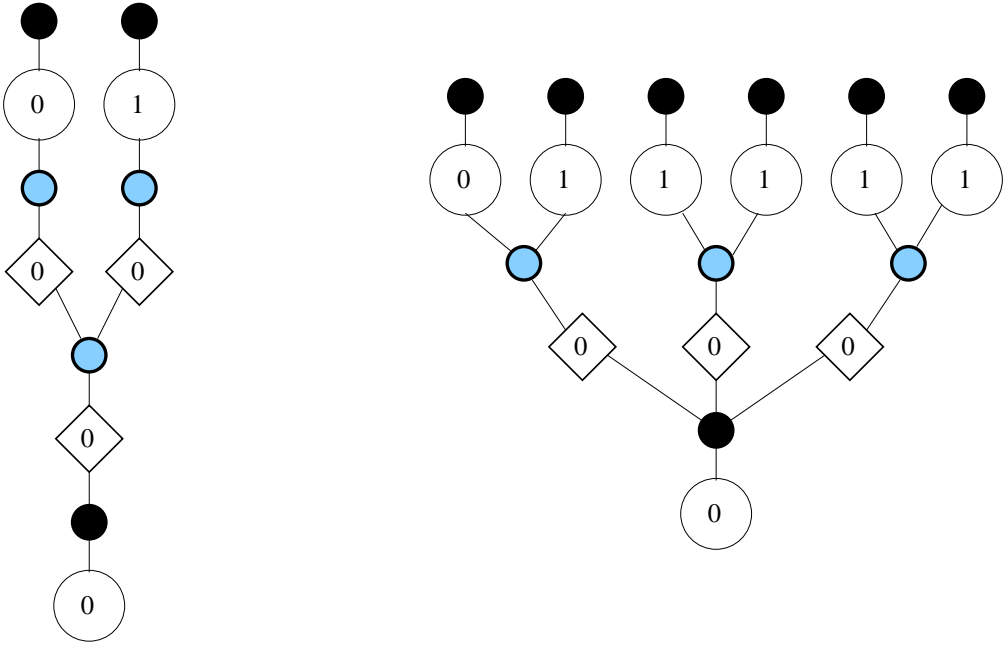}
\caption{The events $F$ (left) and $\Gamma$ for $k=2$ (right)}
\end{figure}

\proof We first treat the case where $q_{2k+1}=0$ for every $k\geq2$
which implies that $q_3>0$. Consider the event \ba
&&E=\Big\{\z_\root(\tr)=(0,0,1,0),\z_1(\tr)=(0,1,0,0),\z_{11}(\tr)=(0,0,0,1),\z_{111}(\tr)=(0,2,0,0),\\
&&\hspace{1.5cm}\z_{1111}(\tr)=\z_{1112}(\tr)=(0,0,0,1),\z_{11111}(\tr)=\z_{11121}(\tr)=(1,0,0,0),\\
&&\hspace{1.5cm}\z_{111111}(\tr)=\z_{111211}(\tr)=(0,0,1,0)\Big\}.
\ea Let $u\in\u$ and let $(\tr,\e,\l)\in\ov{\Tr}_M$ such that
$u\in\tr$. We set $\tr^{[u]}=\{v\in\u:uv\in\tr\}$ and for every
$v\in\tr^{[u]}$ we set $\e^{[u]}(v)=\e(uv)$ and
$\l^{[u]}(v)=\l(uv)-\l(u)$. On the event $E$ we can define
$(\tr_1,\e_1,\l_1)=(\tr^{[u_1]},\e^{[u_1]},\l^{[u_1]})$ and
$(\tr_2,\e_2,\l_2)=(\tr^{[u_2]},\e^{[u_2]},\l^{[u_2]})$, where we have
written $u_1=111111$ and $u_2=111211$. Let $F$ be the event defined by
$$F=E\,\cap\,\Big\{\l_{1}=\l_{11}=\l_{111}=\l_{1111}=\l_{1112}=\l_{11111}=\l_{11121}=\l_{111111}=0,\l_{111211}=1\Big\}.$$
We observe that $\Q_{\mup,\ovlra{\nup}}(F)>0$ and that under
$\Q_{\mup,\ovlra{\nup}}(\cdot\mid F)$, the spatial trees
$(\tr_1,\e_1,\l_1)$ and $(\tr_2,\e_2,\l_2)$ are independant and
distributed according to $\Q_{\mup,\ovlra{\nup}}$. Furthermore
$$\left\{\#\tr^1=n,v_m\in\dd_1\tr\right\}\supset
F\cap\left\{v_{m,1}\in\dd_1\tr_1\right\}\cap\left\{\un{\l}_2\geq0\right\}\cap\left\{\#\tr_1^1+\#\tr_2^1=n-1\right\}\,
,$$ where $v_{m,1}$ is the first vertex of $\tr_1^1\setminus\{\root\}$
that achieves the minimum of $\un{\l}_1$. So we obtain that
\begin{eqnarray}
  &&\hspace{-1cm}\Q_{\mup,\ovlra{\nup}}\left(\#\tr^1=n,v_m\in\dd_1\tr\right)\label{eqargmin1}\\
  &\geq&\Q_{\mup,\ovlra{\nup}}(F)\sum_{j=1}^{n-2}\Q_{\mup,\ovlra{\nup}}\left(\#\tr^1=j,v_m\in\dd_1\tr\right)\Q_{\mup,\ovlra{\nup}}\left(\#\tr^1=n-1-j,\un{\l}\geq0\right).\nonumber
\end{eqnarray}

Let us now turn to the second case for which there exists $k\geq2$
such that $q_{2k+1}>0$. Let $K=2k-1$. On the event \ba
&&\Lambda=\Big\{\z_\root(\tr)=(0,0,1,0),\z_1(\tr)=(0,K,0,0),\z_{11}(\tr)=\ldots=\z_{1K}(\tr)=(0,0,0,1),\\
&&\hspace{1.3cm}\z_{111}(\tr)=\ldots=\z_{1K1}(\tr)=(k,0,0,0),\z_{1111}(\tr)=\z_{1112}(\tr)=(0,0,1,0)\Big\}
\ea we can define
$((\tr^{[u_{ij}]},\e^{[u_{ij}]},\l^{[u_{ij}]}))_{1\leq i\leq K,1\leq
  j\leq k}$ where we have written $u_{ij}=1i1j$. Let $\Gamma$ be the
event
$$\Lambda\cap\{\l_1=0\}\cap\bigcap_{1\leq i\leq
  K}\{\l_{1i}=\l_{1i1}=0\} \cap\{\l_{1111}=0\}\cap\bigcap_{2\leq i\leq
  k}\{\l_{111i}=1\}\cap\hspace{-0.4cm}\bigcap_{2\leq i\leq K,1\leq
  j\leq k}\hspace{-0.4cm}\{\l_{u_{ij}}=1\}.$$ We observe that
$\Q_{\mup,\ovlra{\nup}}(\Gamma)>0$. Furthermore, under the probability
measure $\Q_{\mup,\ovlra{\nup}}(\cdot\mid \Gamma)$, the spatial trees
$((\tr^{[u_{ij}]},\e^{[u_{ij}]},\l^{[u_{ij}]}))_{1\leq i\leq K,1\leq
  j\leq k}$ are independant,
$(\tr^{[u_{11}]},\e^{[u_{11}]},\l^{[u_{11}]})$ and
$(\tr^{[u_{12}]},\e^{[u_{12}]},\l^{[u_{12}]})$ are distributed
according to $\Q_{\mup,\ovlra{\nup}}$, and
$((\tr^{[u_{ij}]},\e^{[u_{ij}]},\l^{[u_{ij}]}))_{1\leq i\leq K,1\leq
  j\leq k}$ are distributed according to
$\P_{\mup,\ovlra{\nup},0}$. Last \ba
\left\{\#\tr^1=n,v_m\in\dd_1\tr\right\}&\supset&
\Gamma\cap\left\{v_m^{u_{11}}\in\dd_1\tr^{[u_{11}]}\right\}\cap
\left\{\un{\l}^{[u_{12}]}\geq0\right\}\\&&\hspace{-1cm}
\cap\left\{\#\tr^{[u_{11}],1}+ \#\tr^{[u_{12}],1}=n+1-kK\right\}
\cap\hspace{-0.2cm}\bigcap_{2\leq i\leq K,1\leq j\leq
  k}\hspace{-0.1cm}\left\{\tr^{[u_{ij}]}=\{\root\}\right\}.  \ea So we obtain that
\begin{eqnarray}
  &&\hspace{-1cm}\Q_{\mup,\ovlra{\nup}}
  \left(\#\tr^1=n,v_m\in\dd_1\tr\right)\label{eqargmin2}\\
  &\geq&C\sum_{j=2}^{n-kK}\Q_{\mup,\ovlra{\nup}}
  \left(\#\tr^1=j,v_m\in\dd_1\tr\right)\Q_{\mup,\ovlra{\nup}}
  \left(\#\tr^1=n+1-kK-j,\un{\l}\geq0\right)\nonumber
\end{eqnarray}
where we have written
$C=\mug^{(1)}(\{\zero\})^{k(K-1)}\Q_{\mup,\ovlra{\nup}}(\Gamma)$.

We can now conclude the proof of Lemma \ref{argmin} in both cases from
respectively (\ref{eqargmin1}) and (\ref{eqargmin2}) by following the
lines of the proof of Lemma 4.3 in \cite{LG}.\cq

We can now state the main result of this section.

\begin{proposition}\label{estpositif} Let $M>0$. There exist four
  constants $\gamma_1>0$, $\gamma_2>0$, $\wt{\gamma}_1>0$ and
  $\wt{\gamma}_2>0$ such that for all $n$ sufficiently large and for
  every $x\in[0,M]$, \ba
  \frac{\wt{\gamma}_1}{n}\;\leq&\Q_{\mup,\ovlra{\nup}}^n(\,\un{\l}>0)&\leq\;\frac{\wt{\gamma}_2}{n},\\
  \frac{\gamma_1}{n}\;\leq&\P_{\mup,\ovlra{\nup},x}^n(\,\un{\l}>0)&\leq\;\frac{\gamma_2}{n}.
  \ea
\end{proposition}

\proof We prove exactly in the same way as in \cite{LG} the existence
of $\wt{\gamma}_2$ and the existence of a constant $\gamma_3>0$ such
that for all $n$ sufficiently large, we have
$$\Q^n_{\mup,\ovlra{\nup}}\left(\un{\l}\geq0\right)\geq\frac{\gamma_3}{n}.$$

Let us now fix $M>0$. Let $k\geq1$ be such that $q_{2k+1}>0$. We
choose an integer $p$ such that $pk\geq M$. First note that
$\ovlra{\nug}_{3,\w}(\{\zero\})=1/(2k-1)$ if $\w=(0,2k-1,0,0)$ and
that $\ovlra{\nug}_{4,\w}(\{k,k-1,\ldots,1\})=1/(2\#A_{k,0})$ if
$\w=(k,0,0,0)$. For every $l\in\N$, we define $1^l\in\u$ by
$1^l=11\ldots1$, $|1^l|=l$. By arguing on the event \ba
E'=\Big\{\z_\root(\tr)=\ldots=\z_{1^{4p}}(\tr)=(0,0,1,0),\z_1(\tr)=\ldots=\z_{1^{4p-3}}(\tr)=(0,2k-1,0,0),\\
\hspace{1cm}\z_{11}(\tr)=\ldots=\z_{1(2k-1)}(\tr)=\ldots=\z_{1^{4p-3}1}(\tr)=\ldots=\z_{1^{4p-3}(2k-1)}=(0,0,0,1),\\
\hspace{1cm}\z_{111}(\tr)=\ldots=\z_{1(2k-1)1}(\tr)=\ldots=\z_{1^{4p-3}11}(\tr)=\ldots=\z_{1^{4p-3}(2k-1)1}=(k,0,0,0)\Big\}\\
\hspace{1.8cm}\cap\bigcap_{i=0}^{p-1}\{\z_{1^{4i+3}2}(\tr)=\ldots=\z_{1^{4i+3}k}=\zero\}\cap\bigcap_{i=0}^{p-1}\bigcap_{j=2}^{2k-1}\bigcap_{l=1}^k\{\z_{1^{4i+1}j1l}=\zero\},
\ea we show that
$$\Q_{\mup,\ovlra{\nup}}\left(\un{\l}>0,\#\tr^1=n\right)\geq\frac{C(\mug,\nug,k)^p}{\mu^{(1)}(\{(0,0,1,0)\})}\,\P_{\mup,\ovlra{\nup},pk}\left(\un{\l}>0,\#\tr^1=n-pk(2k-1)\right),$$
where $C(\mug,\nug,k)$ is equal to
$$\frac{\mug^{(1)}(\{(0,0,1,0)\})\mug^{(3)}(\{(0,2k-1,0,0)\})(\mug^{(4)}(\{(k,0,0,0)\}))^{2k-1}\mug^{(1)}(\{\zero\}))^{k(2k-1)-1}}{(2k-1)(2\#A_{k,0})^{2k-1}}.$$
This implies thanks to Lemma \ref{asympt} that for all $n$ sufficiently large,
$$\P^n_{\mup,\ovlra{\nup},pk}\left(\un{\l}>0\right)\leq\frac{2\mug^{(1)}(\{(0,0,1,0)\})\wt\gamma_2}{C(\mug,\nug,k)n},$$
which ensures the existence of $\gamma_2$.

Last by arguing on the event \ba
&&\hspace{-0.8cm}F=\Big\{\z_\root(\tr)=\z_{1^4}=(0,0,1,0),\z_1(\tr)=(0,2k-1,0,0),\\
&&\z_{11}(\tr)=\ldots=\z_{1(2k-1)}(\tr)=(0,0,0,1),\z_{111}(\tr)=\ldots=\z_{1(2k-1)1}(\tr)=(k,0,0,0)\Big\}\\
&&\hspace{0.8cm}\cap\bigcap_{j=2}^{k}\{\z_{1^3j}(\tr)=\zero\}\cap\bigcap_{i=2}^{2k-1}\bigcap_{j=1}^{k}\{\z_{1i1j}(\tr)=\zero\},
\ea we show that \ba
&&\hspace{-1cm}\P_{\mup,\ovlra{\nup},0}\left(\un{\l}>0,\#\tr^1=n\right)\\
&\geq&
C(\mug,\nug,k)\mug^{(1)}(\{(0,0,1,0)\})\Q_{\mup,\ovlra{\nup}}\left(\un{\l}\geq0,\#\tr^1=n-k(2k-1)\right),
\ea which ensures the existence of $\gamma_1$. We get the existence of
$\wt{\gamma}_1$ by the same arguments.\cq

\subsection{Proofs of Theorem \ref{thconvcond} and of Theorem
  \ref{thcartes}}

To prove Theorem \ref{thconvcond} from what precedes, we can adapt
section 7 of \cite{LG} in exactly the same way as in the proof of
Theorem 3.3 in \cite{We}. A key result in the proof of Theorems 2.2 in
\cite{LG} and 3.3 in \cite{We} is a spatial Markov property for
spatial Galton-Watson trees. Let $a>0$ and $(\tr,\e,\l)\in\Tr_M$. As
in section 5 of \cite{LG} let $v_1,\ldots,v_M$ denote the exit
vertices from $(-\infty,a)$ listed in lexicographical order, and let
$(\tr^a,\e^a,\l^a)$ correspond to the multitype spatial tree
$(\tr,\e,\l)$ which has been truncated at the first exit from
$(-\infty,a)$. Let $v\in\tr$. Recall from section \ref{sec:estpos} the
definition of the multitype spatial tree
$(\tr^{[v]},\e^{[v]},\l^{[v]})$. We set
$\ov{\l}^{[v]}_u=\l^{[v]}_u+\l_v$ for every $u\in\tr^{[v]}$.

\begin{lemma}\label{Markovsp} Let $x\in[0,a)$ and $p\in\{1,\ldots,n\}$. Let $n_1,\ldots,n_p$ be positive integers such that $n_1+\ldots+n_p\leq n$. Assume that
$$\ov{\P}_{\mup,\ovlra\nup,x}^{(1),n}\left(M=p,\;\#\tr^{[v_1],1}=n_1,\ldots,\;\#\tr^{[v_p],1}=n_p\right)>0.$$
Then, under the probability measure $\ov{\P}_{\mup,\nup,x}^{(1),n}(\cdot\mid M=p,\,\#\tr^{[v_1],1}=n_1,\ldots,\,\#\tr^{[v_p],1}=n_p)$, and conditionally on $(\tr^a,\e^a,\l^a)$, the spatial trees 
$$\left(\tr^{[v_1]},\e^{[v_1]},\ov{\l}^{[v_1]}\right),\ldots,\left(\tr^{[v_p]},\e^{[v_p]},\ov{\l}^{[v_p]}\right)$$ 
are independent and distributed respectively according to $\ov{\P}_{\mup,\ovlra\nup,\l_{v_1}}^{(\e(v_1)),n_1},\ldots,\ov{\P}_{\mup,\ovlra\nup,\l_{v_p}}^{(\e(v_p)),n_p}$.
\end{lemma}
Beware that in our context, if $v$ is an exit vertex then $\e(v)\in\{1,2\}$. This is the reason why Theorem \ref{thconvcond} is stated under both probability measures $\P^{(1),n}_{\mup,\ovlra{\nup},x}$ and $\P^{(2),n}_{\mup,\ovlra{\nup},x}$. Thus the statement of Lemma 7.1 of \cite{LG} (and of Lemma 3.18 of \cite{We}) is modified in the following way. Set for every $n\geq1$ and every $s\in[0,1]$,
\ba
C^{(n)}(s)&=&\Ad_\q\frac{C((\#\tr-1)s)}{n^{1/2}},\\
V^{(n)}(s)&=&\Bd_\q\frac{V((\#\tr-1)s)}{n^{1/4}}.
\ea
Last define from section \ref{serpent}, on a suitable probability space $(\OO,\PP)$, a collection of processes $(\ov{\b}^x,\ov{\r}^x)_{x>0}$.

\begin{lemma}\label{condpos}
Let $F:\C([0,1],\R)^2\to\R$ be a Lipschitz function. Let $0<c'<c''$. Then for $i\in\{1,2\}$,
$$\sup_{c'n^{1/4}\leq y\leq c''n^{1/4}}\left|\ov{\E}_{\mup,\nup,y}^{(i),n}\left(F\left(C^{(n)},V^{(n)}\right)\right)-\EE\left(F\left(\ov{\b}^{\Bd_\q y/n^{1/4}},\ov{\r}^{\Bd_\q y/n^{1/4}}\right)\right)\right|\build{\la}_{n\to\infty}^{}0.$$
\end{lemma}

In the remainder of this section we derive Theorem \ref{thcartes} from
Theorem \ref{thconvcond} in the same way as Theorem 2.5 in \cite{We}
is derived from Theorem 3.3. We first state a lemma, which is
analogous to Lemma 3.20 in \cite{We} in our more general setting. To
this end we introduce some notation. Recall that if $\tr\in T_M$,
we set $\xi=\#\tr-1$ and we denote by $v(0)=\root\prec
v(1)\prec\ldots\prec v(\xi)$ the list of the vertices of $\tr$ in
lexicographical order. For $n\in\{0,1,\ldots,\xi\}$, we set as in
\cite{Mi1},
$$\Lambda^\tr_1(n)=\#\left(\tr^1\cap\{v(0),v(1),\ldots,v(n)\}\right).$$
We extend $\Lambda^\tr_1$ to the real interval $[0,\xi]$ by setting
$\Lambda^\tr_1(s)=\Lambda^\tr_1(\lfloor s\rfloor)$ for every
$s\in[0,\xi]$, and we set for every $s\in[0,1]$
$$\ov{\Lambda}^\tr_1(s)=\frac{\Lambda^\tr_1(\xi s)}{\#\tr^1}.$$ 

Recall that $u_0,u_1,\ldots,u_{2\xi}$ denotes the search-depth
sequence of $\tr$. We also define for $k\in\{0,1,\ldots,2\xi\}$,
$$K_\tr(k)=1+\#\left\{l\in\{1,\ldots,k\}:C(l)=C(l-1)+1\;{\rm and}\;\e(u_l)=1\right\}.$$
Note that $K_\tr(k)$ is the number of vertices of type $1$ in the
search-depth sequence up to time $k$. As previously, we extend $K_\tr$
to the real interval $[0,2\xi]$ by setting $K_\tr(s)=K_\tr(\lfloor
s\rfloor)$ for every $s\in[0,2\xi]$, and we set for every $s\in[0,1]$
$$\ov{K}_\tr(s)=\frac{K_\tr(2\xi s)}{\#\tr^1}.$$ 

\begin{lemma}\label{lemISE}
  The law under $\ov{\P}_{\mup,\ovlra{\nup},1}^{(1),n}$ of
  $\left(\ov{\Lambda}^\tr_1(s),0\leq s\leq1\right)$ converges as
  $n\to\infty$ to the Dirac mass at the identity mapping of
  $[0,1]$. In other words, for every $\eta>0$, \be\label{eqISEMi1}
  \ov{\P}_{\mup,\ovlra\nup,1}^{(1),n}\left(\sup_{s\in[0,1]}\left|\ov{\Lambda}^\tr_1(s)-s\right|>\eta\right)\build{\la}_{n\to\infty}^{}0.
  \ee Consequently, the law under
  $\ov{\P}_{\mup,\ovlra\nup,1}^{(1),n}$ of $\left(\ov{K}_\tr(s),0\leq
    s\leq1\right)$ converges as $n\to\infty$ to the Dirac mass at the
  identity mapping of $[0,1]$. In other words, for every $\eta>0$,
  \be\label{eqISEK}
  \ov{\P}_{\mup,\ovlra\nup,1}^{(1),n}\left(\sup_{s\in[0,1]}\left|\ov{K}_\tr(s)-s\right|>\eta\right)\build{\la}_{n\to\infty}^{}0.
  \ee
\end{lemma}

\proof For $\tr\in T_M$, we let $v^1(0)=\root\prec
v^1(1)\prec\ldots\prec v^1(\#\tr^1-1)$ be the list of vertices of
$\tr$ of type $1$ in lexicographical order. We define as in \cite{Mi1}
$$G^\tr_1(k)=\#\left\{u\in\tr:u\prec v^1(k)\right\},\;\;0\leq k\leq\#\tr^1-1,$$
and we set $G^\tr_1(\#\tr^1)=\#\tr$. Note that $v^1(k)$ does not
belong to the set $\{u\in\tr:u\prec v^1(k)\}$. Recall from section \ref{secdeflois} the definition of
the vector $\ag=(a_1,a_2,a_3,a_4)$. From the second assertion of
Proposition 6 in \cite{Mi1}, for every $s\in[0,1]$, there exists a
constant $\vep>0$ such that for all $n$ sufficiently large,
$$P_\mup^{(1)}\left(|G^\tr_1(\lfloor ns\rfloor)-a_1^{-1}ns|\geq n^{3/4}\right)\leq e^{-n^\vep}.$$
Thus we obtain thanks to Lemma \ref{asympt} and Proposition
\ref{estpositif} that for every $s\in[0,1]$, there exists a constant
$\vep'>0$ such that for all $n$ sufficiently large,
$$\ov\P_{\mup,\ovlra\nup,1}^{(1),n}\left(|G^\tr_1(\lfloor ns\rfloor)-a_1^{-1}ns|\geq n^{3/4}\right)\leq e^{-n^\vep}.$$
Let us fix $\eta>0$. We then have for every $s\in[0,1]$,
$$\ov\P_{\mup,\ovlra\nup,1}^{(1),n}\left(|n^{-1}G^\tr_1(\lfloor ns\rfloor)-a_1^{-1}s|\geq\eta\right)\build{\la}_{n\to\infty}^{}0.$$
In particular for $s=1$ we have
$$\ov\P_{\mup,\ovlra\nup,1}^{(1),n}\left(|n^{-1}\#\tr-a_1^{-1}|\geq\eta\right)\build{\la}_{n\to\infty}^{}0,$$
which implies that for every $s\in[0,1]$,
$$\ov\P_{\mup,\ovlra\nup,1}^{(1),n}\left(|(\#\tr)^{-1}G^\tr_1(\lfloor ns\rfloor)-s|\geq\eta\right)\build{\la}_{n\to\infty}^{}0.$$
Let us now set $k_\eta=\lceil\frac{2}{a_1\eta}\rceil$ and
$s_m=mk_\eta^{-1}$ for every $m\in\{0,1,\ldots,k_\eta\}$. Since the
mapping $s\in[0,1]\mapsto n^{-1}G^\tr_1(\lfloor ns\rfloor)$ is
non-decreasing, we have
$$\ov\P_{\mup,\ovlra\nup,1}^{(1),n}\left(\sup_{s\in[0,1]}\left|\frac{G^\tr_1(\lfloor ns\rfloor)}{\#\tr}-s\right|\geq\eta\right)\leq\ov\P_{\mup,\ovlra\nup,1}^{(1),n}\left(\sup_{0\leq m\leq k_\eta}\left|\frac{G^\tr_1(\lfloor ns_m\rfloor)}{\#\tr}-s_m\right|\geq\frac{\eta}{2}\right),$$
implying that
$$\ov\P_{\mup,\ovlra\nup,1}^{(1),n}\left(\sup_{s\in[0,1]}\left|(\#\tr)^{-1}G^\tr_1(\lfloor ns\rfloor)-s\right|\geq\eta\right)\build{\la}_{n\to\infty}^{}0.$$
We thus get (\ref{eqISEK}) in the same way as (32) is obtained in
\cite{We}. Then we derive (\ref{eqISEK}) from (\ref{eqISEMi1}) in the
same way as (33) is derived from (32) in \cite{We}.\cq

We can now complete the proof of Theorem \ref{thcartes}. Recall that
$\rad_\map$ denotes the radius of the map $\map$. Thanks to
Proposition \ref{loiimage} we know that the law of $\rad_\map$ under
$\B_\q^r(\cdot\mid\#\v_\map=n)$ coincides with the law of
$\sup_{v\in\tr^1}\l_v$ under $\ov\P_{\mup,\nup,1}^{(1),n}$. Furthermore
we easily see (compare \cite[Lemma 1]{Mi2}) that the law of
$\sup_{v\in\tr^1}\l_v$ under $\ov\P_{\mup,\nup,1}^{(1),n}$ is the law of
$\sup_{v\in\tr^1}\l_v$ under $\ov\P_{\mup,\ovlra\nup,1}^{(1),n}$. We thus
get the first assertion of Theorem \ref{thcartes}.

Let us turn to (ii). By Proposition \ref{loiimage} and properties of
the Bouttier-Di Francesco-Guitter bijection, the law of
$\lm_\map^{(n)}$ under $\B_\q^r(\cdot\mid\#\v_\map=n)$ is the law
under $\ov{\P}_{\mup,\nup,1}^{(1),n}$ of the probability measure
$\I_n$ defined by
$$\langle\I_n,g\rangle=\frac{1}{\#\tr^1+1}
\left(g(0)+\sum_{v\in\tr^1}g\left(n^{-1/4}\l_v\right)\right),$$ which
coincides with the law of $\I_n$ under
$\ov\P_{\mup,\ovlra\nup,1}^{(1),n}$. We thus complete the proof of (ii)
by following the lines of the proof of Theorem 2.5 in \cite{We}. Last assertion (iii) easily follows from (ii).


\end{document}